\documentclass[12pt,english,twoside]{article}
\usepackage[english]{babel}
\usepackage{amsmath,amsfonts}
\usepackage{mathtools}
\usepackage[cp850]{inputenc}
\usepackage{tikz}
\allowdisplaybreaks
\newcounter{theorem}[section]
\renewcommand{\thetheorem}{\thesection.\arabic{theorem}}
\newenvironment{lemma}[      1]{\refstepcounter{theorem} %
\bf \thetheorem\ Lemma#1.        \it}{}
\newenvironment{theorem}[    1]{\refstepcounter{theorem} %
\bf \thetheorem\ Theorem#1.      \it}{}

\newenvironment{corollary}[  1]{\refstepcounter{theorem} %
\bf \thetheorem\ Corollary#1.    \it}{}
\newenvironment{example}[    1]{\refstepcounter{theorem} %
\bf \thetheorem\ Example#1.      \rm}{\par}
\newenvironment{examples}[   1]{\refstepcounter{theorem} %
\bf \thetheorem\ Examples#1.     \rm}{\par}

\newenvironment{proof}{ %
\it              Proof.          \rm}{\hfill $ \Box $}

\newcounter{abc}[theorem]

\newcounter{one}[theorem]
\newenvironment{onelist}{\begin{list}{%
\rm (\arabic{one}) \hfill           }{\usecounter{one} %
\topsep0mm \partopsep0mm \parsep0mm \itemsep0mm %
\leftmargin2em \labelwidth2em \labelsep0em}}{\end{list}}

\newenvironment{hylist}{\begin{list}{%
--               \hfill            }{%
\topsep0mm \partopsep0mm \parsep0mm \itemsep0mm %
\leftmargin2em \labelwidth2em \labelsep0em}}{\end{list}}

\newcounter{rom}
\newenvironment{romlist}{\begin{list}{%
\rm (\roman{rom})  \hfill           }{\usecounter{rom} %
\topsep0mm \partopsep0mm \parsep0mm \itemsep0mm %
\leftmargin2em \labelwidth2em \labelsep0em}}{\end{list}}

\newcommand{\fa}[3]{\{#1_#2\}_{#2\in#3}}

\newcommand{\I}{{\mathbb I}}
\newcommand{\N}{{\mathbb N}}
\newcommand{\R}{{\mathbb R}}

\newcommand{\BB}{{\cal B}}
\newcommand{\CC}{{\cal C}}

\newcommand{\aaa}{{\bf a}}
\newcommand{\bbb}{{\bf b}}

\newcommand{\eee}{{\bf e}}

\newcommand{\sss}{{\bf s}}

\newcommand{\uuu}{{\bf u}}
\newcommand{\vvv}{{\bf v}}
\newcommand{\www}{{\bf w}}
\newcommand{\xxx}{{\bf x}}
\newcommand{\yyy}{{\bf y}}
\newcommand{\zzz}{{\bf z}}

\newcommand{\FFF}{{\bf F}}

\newcommand{\HHH}{{\bf H}}

\newcommand{\XXX}{{\bf X}}

\newcommand{\dps}{\displaystyle}

\newcommand{\leb}{{\lambda}}
\newcommand{\eeta}{{\boldsymbol\eta}}

\newcommand{\zero}{{\bf 0}}
\newcommand{\eins}{{\bf 1}}
\newcommand{\iinfty}{{\boldsymbol \infty}}

\newcommand{\leftmatrix}[1]{\left(\begin{array}{#1}}
\newcommand{\rightmatrix}{\end{array}\right)}

\DeclareMathOperator*{\bcirc}{\bigcirc}

\begin{document}

\title{\Large\bf On Kendall's Tau for Order Statistics}
\author{\normalsize Sebastian Fuchs\footnote{E--mail Sebastian Fuchs: {\tt sfuchs@statistik.tu-dortmund.de} (Corresponding author)} \\[1ex]
\small Fakult{\"a}t Statistik \\[-0.8ex]
\small Technische Universit{\"a}t Dortmund \\[-0.8ex] 
\small 44221 Dortmund \\[-0.8ex] 
\small Germany 
\and 
\normalsize Klaus D. Schmidt\footnote{E--mail Klaus D. Schmidt: {\tt klaus.d.schmidt@uni-mannheim.de}} \\[1ex]
\small Institut f{\"u}r Mathematik \\[-0.8ex]
\small Universit{\"a}t Mannheim \\[-0.8ex]
\small 68161 Mannheim \\[-0.8ex] 
\small Germany }
\date{}
\maketitle

\begin{abstract}
\noindent
Every copula   $ C $   for a random vector   $ \XXX=(X_1,\dots,X_d) $   with identically distributed coordinates 
determines a unique copula   $ C_{:d} $   for its order statistic   $ \XXX_{:d}=(X_{1:d},\dots,X_{d:d}) $. 
In the present paper we study the dependence structure of   $ C_{:d} $   via Kendall's tau, 
denoted by   $ \kappa $. 
As a general result, 
we show that   $ \kappa[C_{:d}] $   is at least as large as   $ \kappa[C] $. 
For the product copula   $ \Pi $, 
which corresponds to the case of independent coordinates of   $ \XXX $, 
we provide an explicit formula for   $ \kappa[\Pi_{:d}] $   showing that the inequality between   $ \kappa[\Pi] $   and   $ \kappa[\Pi_{:d}] $   is strict. 
We also compute Kendall's tau for certain multivariate margins of   $ \Pi_{:d} $   corresponding to the lower or upper coordinates of   $ \XXX_{:d} $. 
\end{abstract}



\newpage


\section{Introduction}
\label{intro}

Kendall's tau is a measure of association which evaluates every copula   $ C : [0,1]^d\to[0,1] $   by a real number   $ \kappa[C] $   
satisfying   $ - 1/(2^{d-1}\!-\!1) \leq \kappa[C] \leq 1 $. 

\bigskip
In the present paper we study Kendall's tau for a class of copulas related to the order statistic   $ \XXX_{:d} = (X_{1:d},\dots,X_{d:d}) $   
of a random vector   $ \XXX = (X_1,\dots,X_d) $   with identical univariate marginal distributions.  
In this case, 
every copula   $ C $   for   $ \XXX $   determines a copula   $ C_{:d} $   for   $ \XXX_{:d} $; 
see 
Navarro and Spizzichino [2010] for the special case where the distribution functions of the coordinates of   $ \XXX $   are continuous and 
Dietz et al.\ [2016] for the general case. 
Since the construction of the order transform   $ C_{:d} $   of a copula   $ C $   
is determined by a map transforming every random vector into its order statistic, 
we shall study Kendall's tau for the order transform   $ C_{:d} $   of an arbitrary copula   $ C $. 

\bigskip
As a general result, 
we show that the order transform   $ C_{:d} $   of a copula   $ C $   satisfies   $ \kappa[C] \leq \kappa[C_{:d}] $   
(Theorem 
\ref{t.Kendall}). 
We also show that the inequality is strict when   $ C $   is the product copula 
(Theorem 
\ref{t.product}), 
which corresponds to the case where the coordinates of   $ \XXX $   are also independent. 
By contrast, 
the inequality becomes an equality when   $ C $   is the upper Fr{\'e}chet--Hoeffding bound (and hence maximizes Kendall's tau) 
or a copula which is symmetric and minimizes Kendall's tau 
(Corollary 
\ref{c.bounds} and Theorem 
\ref{t.Kendall}). 
This result is of interest 
since the computation of   $ C_{:d} $   or   $ \kappa[C_{:d}] $   may be tedious and is not needed in the case where   $ \kappa[C] = \kappa[C_{:d}] $   and 
since there exist many symmetric copulas minimizing Kendall's tau; 
see 
Fuchs et al.\ [2018]. 

\bigskip
The major part of this paper is devoted to the order transform   $ \Pi_{:d} $   of the product copula   $ \Pi $.   
We first determine   $ \kappa[\Pi_{:d}] $   
(Theorem 
\ref{t.product}) and then present some identities for the computation of   $ \kappa[\varrho_K(\Pi_{:d})] $,   
where   $ \varrho_K(\Pi_{:d}) $   denotes the multivariate margin of   $ \Pi_{:d} $   with respect to the coordinates 
in   $ K\subseteq\{1,\dots,d\} $   with   $ |K|\geq2 $   
(Theorem 
\ref{t.product-K}). 
In particular, 
we obtain an explicit formula for Kendall's tau of   $ \varrho_{\{1,\dots,k\}}(\Pi_{:d}) $   
(Corollary 
\ref{c.product-lowertail}), 
which is a copula for the lower   $ k $   coordinates of   $ \XXX_{:d} $, 
and we show that, 
due to a general reflection principle 
(Theorem 
\ref{t.product-reflection}), 
this formula is also valid for Kendall's tau of   $ \varrho_{\{d-k+1,\dots,d\}}(\Pi_{:d}) $,   
which is a copula for the upper   $ k $   coordinates of   $ \XXX_{:d} $. 
We thus extend certain results for the case   $ |K|=2 $; 
see 
Av{\'e}rous et al.\ [2005] and 
Navarro and Balakrischnan [2010]. 

\bigskip
This paper is organized as follows: 
Section \ref{preliminaries} collects some definitions and results on copulas and related topics which will be needed in this paper. 
Section \ref{kendall} provides a brief discussion of Kendall's tau and Kendall's distribution function. 
Section \ref{ot-copula} starts with the definition of the order transform of the Euclidean space, 
which turns every random vector into its order statistic, 
and proceeds with the construction of the order transform of a copula.  
In Section \ref{ot-kendall} we present some result of Kendall's tau for the order transform of an arbitrary copula and 
in Section \ref{ot-product} we study Kendall's for the order transform of the product copula and its multivariate margins. 
Some auxiliary results needed in Section \ref{ot-product} are established in the Appendix. 

\bigskip
Throughout this paper we shall use the following notation: 
Let   $ \I := [0,1] $   and 
let   $ \leb $   denote the Lebesgue measure on   $ \BB(\R) $. 
Furthermore, 
let   $ d\geq2 $   be an integer, 
which will be kept fixed, 
and 
let   $ \leb^d $   denote the Lebesgue measure on   $ \BB(\R^d) $. 
We denote 
by   $ \eee_1,\dots,\eee_d $   the standard basic unit vectors in   $ \R^d $, 
by   $ \zero $   the vector in   $ \R^d $   with all coordinates being equal to   $ 0 $   and 
by   $ \eins $   the vector in   $ \R^d $   with all coordinates being equal to   $ 1 $.
For   $ \xxx,\yyy\in\R^d $, 
we write   $ \xxx\leq\yyy $   if   $ x_k \leq y_k $   holds for every   $ k\in\{1,\dots,d\} $. 
Then we have   $ \I^d = [\zero,\eins] $. 
On the collection of all real--valued maps on a set   $ S $, 
the pointwise order   $ \leq $   is defined by letting   $ g \leq h $   if   $ g(s) \leq h(s) $   holds for every   $ s \in S $. 

\bigskip
Due to the central role of the order transform   $ T $   of the Euclidean space in the construction of the order transform of a copula   $ C $, 
the symbol   $ C_{:d} $   used in the Abstract and in this Introduction will henceforth be replaced by   $ C_T $.


\section{Preliminaries}
\label{preliminaries}

In this section, 
we recall some definitions and results on 
copulas, 
copula measures, 
groups of transformations of copulas and 
a biconvex form for copulas. 
For further details we refer to 
Fuchs [2014; 2016].

\subsection*{Copulas}

For   $ K\subseteq\{1,...,d\} $, 
we consider the map   $ \eeta_K : \I^d\times\I^d \to \I^d $   given coordinatewise by
\begin{eqnarray*}
        (\eeta_K(\uuu,\vvv))_k                
& := &  \begin{cases}                         
         u_k & k \in \{1,...,d\} \setminus K  \\
         v_k & k \in K                        
        \end{cases}                           
\end{eqnarray*}
and for   $ k\in\{1,\dots,d\} $   we put   $ \eeta_k := \eeta_{\{k\}} $. 

\bigskip
A 
\emph{copula} is a function   $ C: \I^{d} \to \I $   satisfying the following conditions: 
\begin{romlist}
\item   The inequality   $ \sum_{K \subseteq \{1,...,d\}} (-1)^{d-|K|}\,C(\eeta_K(\uuu,\vvv)) \geq 0 $   
holds for all   $ \uuu,\vvv\in\I^d $   such that   $ \uuu\leq\vvv $. 
\item   The identity   $ C(\eeta_k(\uuu,\zero)) = 0   $   holds for every   $ k\in\{1,...,d\} $   and every   $ \uuu\in\I^d $. 
\item   The identity   $ C(\eeta_k(\eins,\uuu)) = u_k $   holds for every   $ k\in\{1,...,d\} $   and every   $ \uuu\in\I^d $. 
\end{romlist}
This definition of a copula is in accordance with the literature;
see, 
e.\,g., 
Durante and Sempi [2016] and 
Nelsen [2006]. 
The collection   $ \CC $   of all copulas is convex.

\bigskip
The following copulas are of particular interest: 
\begin{hylist}
\item   The 
\emph{upper Fr{\'e}chet--Hoeffding bound}   $ M $   given by 
$ M(\uuu) := \min\{u_1,\dots,u_d\} $ 
is a copula and every copula   $ C $   satisfies   $ C \leq M $. 
\item   The 
\emph{product copula}   $ \Pi $   given by 
$ \Pi(\uuu) := \prod_{k=1}^d u_k $ 
is a copula. 
\item   In the case   $ d=2 $, 
the 
\emph{lower Fr{\'e}chet--Hoeffding bound}   $ W $   given by 
$ W(\uuu) := \max\{ \sum_{k=1}^d u_k + 1 - d, 0 \} $   
is a copula and every copula   $ C $   satisfies   $ W \leq C $. 
\end{hylist}

\subsection*{A Group of Transformations of Copulas}

Let   $ \Phi $   denote the collection of all transformations   $ \CC\to\CC $   
and consider the composition   $ \circ: \Phi\times\Phi \to \Phi $   given by   $ (\varphi_2\circ\varphi_1)(C) := \varphi_2(\varphi_1(C)) $   
and the map   $ \iota\in\Phi $   given by   $ \iota(C) := C $. 
Then   $ (\Phi,\circ) $   is a semigroup with neutral element   $ \iota $. 
For   $ i,j,k\in\{1,...,d\} $   with   $ i \neq j $, 
we define the maps   $ \pi_{i,j},\nu_k : \CC\to\CC $   by letting
\begin{eqnarray*}
        (\pi_{i,j}(C))(\uuu)                                      
& := &  C(\eeta_{\{i,j\}}(\uuu,u_j\,\eee_i+u_i\,\eee_j))          \\*
        (\nu_k    (C))(\uuu)                                      
& := &  C(\eeta_k(\uuu,\eins)) - C(\eeta_k(\uuu,\eins\!-\!\uuu))  
\end{eqnarray*}
Each of these maps is an involution and there exists 
\begin{hylist}
\item   a smallest subgroup   $ \Gamma^\pi $   of   $ \Phi $   containing every   $ \pi_{i,j} $, 
\item   a smallest subgroup   $ \Gamma^\nu $   of   $ \Phi $   containing every   $ \nu_k $   and 
\item   a smallest subgroup   $ \Gamma     $   of   $ \Phi $   containing   $ \Gamma^\pi\cup\Gamma^\nu $. 
\end{hylist}
The group   $ \Gamma^\nu $   is commutative. 
Moreover, 
the 
\emph{total reflection} 
\begin{eqnarray*}
        \tau                              
& := &  \bcirc_{k\in\{1,...,d\}} \nu_{k}  
\end{eqnarray*}
transforms every copula into its survival copula and satisfies   $ \tau(M)=M $, 
and we put   $ \Gamma^\tau := \{\iota,\tau\} $. 
The total reflection is used in the definition of the 
\emph{concordance order}   $ \leq_c $   on   $ \CC $   
which is defined by letting   $ C \leq_c D $  if and only if 
$ C \leq D $   and   $ \tau(C) \leq \tau(D) $. 

\bigskip
The group   $ \Gamma $    is a representation of the hyperoctahedral group, 
which also has a well--known geometric representation: 

\bigskip
Let   $ \tilde\Phi $   denote the collection of all transformations   $ \I^d\to\I^d $   
and consider the composition   $ \diamond: \I^d\times\I^d\to\I^d $   
given by   $ (\tilde\varphi_2\diamond\tilde\varphi_1)(\uuu) := \tilde\varphi_2(\tilde\varphi_1(\uuu)) $   
and the map   $ \tilde\iota\in\Phi $   
given by   $ \tilde\iota(\uuu) := \uuu $. 
Then   $ (\tilde\Phi,\diamond) $ is a semigroup with neutral element   $ \tilde\iota $. 
For   $ i,j,k\in\{1,...,d\} $   with   $ i \neq j $, 
we define the maps   $ \tilde\pi_{i,j},\tilde\nu_k : \I^d\to\I^d $   by letting
\begin{eqnarray*}
        \tilde\pi_{i,j}(\uuu)                           
& := &  \eeta_{\{i,j\}}(\uuu,u_j\,\eee_i+u_i\,\eee_j))  \\*
        \tilde\nu_k    (\uuu)                           
& := &  \eeta_k(\uuu,\eins\!-\!\uuu)                    
\end{eqnarray*}
Each of these maps is an involution and there exists 
\begin{hylist}
\item   a smallest subgroup   $ \tilde\Gamma^\pi $   of   $ \tilde\Phi $   containing every   $ \tilde\pi_{i,j} $, 
\item   a smallest subgroup   $ \tilde\Gamma^\nu $   of   $ \tilde\Phi $   containing every   $ \tilde\nu_k $   and 
\item   a smallest subgroup   $ \tilde\Gamma     $   of   $ \tilde\Phi $   containing   $ \tilde\Gamma^\pi\cup\tilde\Gamma^\nu $. 
\end{hylist}
The groups   $ \Gamma $   and   $ \tilde\Gamma $   are related by an isomorphism   $ J : (\Gamma,\circ) \to (\tilde\Gamma,\diamond) $   
satisfying   $ J(\pi_{i,j}) = \tilde\pi_{i,j} $   and   $ J(\nu_k) = \tilde\nu_k $   for all   $ i,j,k\in\{1,...,d\} $. 
For   $ \gamma\in\Gamma $   we put   $ \tilde\gamma := J(\gamma) $, 
and for   $ \tilde\gamma\in\tilde\Gamma $   we put   $ \gamma := J^{-1}(\tilde\gamma) $. 
Then   $ \pi(C) = C\circ\tilde\pi $   holds for every   $ \pi\in\Gamma^\pi $   and every   $ C\in\CC $.

\subsection*{Copula Measures}

Since every copula   $ C\in\CC $   has a unique extension to a distribution function   $ \R^d\to\I $, 
there exists a unique probability measure   $ Q^C : \BB(\I^d)\to\I $   satisfying   
\begin{eqnarray*}
        Q^C[[\zero,\uuu]]  
&  = &  C(\uuu)            
\end{eqnarray*}
for every   $ \uuu\in\I^d $. 
The probability measure   $ Q^C $   is said to be the 
\emph{copula measure} with respect to   $ C $. 
It satisfies   $ Q^C[(\uuu,\vvv)] = Q^C[[\uuu,\vvv]] $   for all   $ \uuu,\vvv\in\I^d $   such that   $ \uuu\leq\vvv $, 
and   $ Q^{\gamma(C)} = (Q^C)_{\tilde\gamma} $   holds for every   $ \gamma\in\Gamma $.

\subsection*{A Biconvex Form for Copulas}

Consider the map   $ [.\,,.]: \CC\times\CC \to \R $   given by 
\begin{eqnarray*}
        [C,D]                            
& := &  \int_{\I^d} C(\uuu)\,dQ^D(\uuu)  
\end{eqnarray*}
The map   $ [.\,,.] $   is in either argument linear with respect to convex combinations and is therefore called a 
\emph{biconvex form}. 
Moreover, 
the map   $ [.\,,.] $   is in either argument monotone with respect to the concordance order. 
It satisfies   $ 0 \leq [C,D] \leq 1/2 $   for all   $ C,D\in\CC $, 
and the bounds are attained since   $ [M,M] = 1/2 $   and since   $ [\nu(M),\nu(M)] = 0 $   holds for every   $ \nu\in\Gamma^\nu\setminus\Gamma^\tau $. 
We also note that   $ [\tau(C),\tau(D)] = [D,C] $   holds for all   $ C,D\in\CC $.


\section{Kendall's Tau and Kendall's Distribution Function}
\label{kendall}

The map   $ \kappa : \CC\to\R $   given by 
\begin{eqnarray*}
        \kappa[C]                                                       
& := &  \frac{2^{d}}{2^{d-1}-1}\,\biggl( [C,C] - \frac{1}{2^d} \biggr)  
\;\,=\,\;  \frac{2^d\,[C,C]-1}{2^{d-1}-1}
\end{eqnarray*}
is called 
\emph{Kendall's tau}; 
see 
Nelsen [2002]. 
Kendall's tau satisfies 
\begin{eqnarray*}
        -\,\frac{1}{2^{d-1}-1}  
&\leq&  \kappa[C]               
\;\,\leq\,\;  1                 
\end{eqnarray*}
and the bounds are attained. 
The following result is obvious from the properties of the biconvex form and will be tacitly used throughout this paper: 

\bigskip
\begin{lemma}{}
\label{l.maxmin}
\begin{onelist}
\item   Kendall's tau is monotone with respect to the concordance order. 
\item   A copula   $ C $   maximizes Kendall's tau if and only if   $ [C,C] = 1/2 $. 
\item   A copula   $ C $   minimizes Kendall's tau if and only if   $ [C,C] = 0   $. 
\end{onelist}
\end{lemma}

\bigskip
For a copula   $ C $, 
the function   $ K_C: \R\to\I $   given by 
\begin{eqnarray*}
        K_C(t)                                        
& := &  Q^C[\{ \uuu\in\I^d \colon C(\uuu) \leq t \}]  
\;\,=\,\;  (Q^C)_C[[0,t]]                             
\end{eqnarray*}
is said to be 
\emph{Kendall's distribution function} with respect to   $ C $, 
which is well--known from the literature. 
It is easy to see that   $ K_C $   is indeed a distribution function 
which satisfies   $ K_C(0) = 0 $   and   $ K_C(1) = 1 $   as well as   $ t \leq  K_C(t) $   for every   $ t\in\I $. 
The following result implies that Kendall's tau can be expressed in term of Kendall's distribution function: 

\bigskip
\begin{lemma}{}
\label{distribution.l-Kendall-2}
The identity 
$$  [C,C]  =  \int_{\I} \Bigl( 1 - K_C(t) \Bigr) d\leb(t)  $$
holds for every copula   $ C $. 
In particular, 
$ \kappa[C] \leq \kappa[D] $   holds for any two copulas   $ C $   and   $ D $   satisfying   $ K_D \leq K_C $. 
\end{lemma}

\bigskip
We refer to Fuchs et al.\ [2018] for further details on the results of this section. 
The comparison of copulas via Kendall's distribution function was considered by 
Cap{\'e}ra{\`a} et al.\ [1997].


\section{The Order Transform of a Copula} 
\label{ot-copula}

Throughout this section, 
consider the map   $ T : \overline\R^d\to\overline\R^d $   which is defined coordinatewise by letting 
\begin{eqnarray*}
        (T(\xxx))_k                                                     
& := &  \min_{J\subseteq\{1,\dots,d\},\;|J|=k} \; \max_{l \in J}\; x_l  
\end{eqnarray*}
Then the coordinates of   $ T(\xxx) $   satisfy   $ (T(\xxx))_1 \leq \dots \leq (T(\xxx))_d $, 
and  $ \xxx\leq\yyy $   implies   $ T(\xxx) \leq T(\yyy) $. 
Moreover, 
the map   $ T $   is measurable and satisfies   $ T^{-1}(\I^d) = \I^d $. 
The map   $ T $   is called the 
\emph{order transform} and associates with every random vector its 
\emph{order statistic}. 

\bigskip
For a distribution function   $ F : \overline\R^d\to\I $, 
let   $ F_1,\dots,F_d $   denote the univariate marginal distribution functions of   $ F $   and 
let   $ F_1^\leftarrow,\dots,F_d^\leftarrow $   denote the corresponding lower quantile functions. 
Consider also the maps   $ \FFF : \overline\R^d\to\I^d $   and   $ \FFF^\leftarrow : \I^d\to\overline\R^d $   which are defined coordinatewise 
by letting 
$$   (\FFF(\xxx))_k             
 :=  F_k(x_k)                   
\qquad\text{and}\qquad          
     (\FFF^\leftarrow(\uuu))_k  
 :=  F_k^\leftarrow(u_k)        
$$
Furthermore, 
let   $ Q^F $   denote the distribution   $ \BB(\overline\R^d) \to \I $   corresponding to   $ F $   and 
let   $ F_T $   denote the distribution function corresponding to   $ (Q^F)_T $. 
Then we have   $ Q^{F_T} = (Q^F)_T $. 

\bigskip
For the remainder of this section, 
consider a fixed copula   $ C $   and let   $ H^C $   denote the distribution function extending   $ C $. 
Then we have 
\begin{eqnarray*}
        H^C           
&  = &  C\circ\HHH^C  
\end{eqnarray*}
Since the univariate marginal distribution functions of   $ H^C   $   are continuous, 
those of   $ H^C_T := (H^C)_T $   are continuous as well. 
Therefore, 
there exists a unique copula   $ C_T $   satisfying 
\begin{eqnarray*}
        H^C_T             
&  = &  C_T\circ\HHH^C_T  
\end{eqnarray*}
and hence   $ C_T = H^C_T \circ (\HHH^C_T)^\leftarrow $. 
The copula   $ C_T $   was introduced and studied by 
Dietz et al.\ [2016] and is called the 
\emph{order transform} of the copula   $ C $. 

\bigskip
The following theorem is due to 
Dietz et al.\ [2016; Theorem 5.2]: 

\bigskip
\begin{theorem}{}
\label{equal.t}
Let   $ F $   be a distribution function satisfying   $ F_1=\ldots=F_n $.   
If   $ C $   is a copula for   $ F $, 
then   $ C_T $   is a copula for   $ F_T $. 
\end{theorem}

\bigskip
We complete this section with two technical results which will be needed later. 

\bigskip
\begin{lemma}{}
\label{l-image}
$ H^C_T $   satisfies   
$ (Q^{H^C_T\circ(\HHH^C_T)^\leftarrow})_{(\HHH^C_T)^\leftarrow} = Q^{H^C_T} $   and 
$ (Q^{C_T})_{C_T} = (Q^{H^C_T})_{H^C_T} $. 
\end{lemma}%

\bigskip
\begin{proof}
For every   $ \xxx\in\overline\R^d $   we have 
\begin{eqnarray*}
        (Q^{H^C_T\circ(\HHH^C_T)^\leftarrow})_{(\HHH^C_T)^\leftarrow}[[-\iinfty,\xxx]]                                                         
&  = &   Q^{H^C_T\circ(\HHH^C_T)^\leftarrow} \Bigl[ \Bigl\{ \uuu\in\I^d \Bigm| (\HHH^C_T)^{\leftarrow}(\uuu)\leq         \xxx  \Bigr\} \Bigr]  \\
&  = &   Q^{H^C_T\circ(\HHH^C_T)^\leftarrow} \Bigl[ \Bigl\{ \uuu\in\I^d \Bigm|                         \uuu \leq\HHH^C_T(\xxx) \Bigr\} \Bigr]  \\[.5ex]
&  = &  (   H^C_T\circ(\HHH^C_T)^\leftarrow \circ \HHH^C_T )(\xxx)                                                                             \\[1ex]
&  = &      H^C_T                                           (\xxx)                                                                             \\*[1ex]
&  = &   Q^{H^C_T}                                                   [[-\iinfty,\xxx]]                                                         
\end{eqnarray*}
This yields the first identity. 
Since   $ C_T = H^C_T\circ(\HHH^C_T)^\leftarrow $, 
the second identity follows from the first. 
\end{proof}

\bigskip
\begin{lemma}{}
\label{l.commute}
Let   $ S : \overline\R^d\to\overline\R^d $   be a measurable map satisfying   $ S(\I^d)\subseteq\I^d $   and 
$$  ((\HHH^C_T)^\leftarrow \circ S)(\uuu)  =  (S \circ (\HHH^C_T)^\leftarrow)(\uuu)  $$
for every   $ \uuu\in\I^d $. 
Then 
$$  \int_{\I^d}          (  C_T \circ S)(\uuu)\,dQ^{  C_T}(\uuu)  
 =  \int_{\overline\R^d} (H^C_T \circ S)(\xxx)\,dQ^{H^C_T}(\xxx)  
$$%
\end{lemma}%

\bigskip
\begin{proof}
Lemma 
\ref{l-image} yields 
\begin{eqnarray*}
        \int_{\I^d}          (  C_T                             \circ S)(\uuu)\,d Q^{  C_T}                          (\uuu)                           
&  = &  \int_{\I^d}          (H^C_T \circ (\HHH^C_T)^\leftarrow \circ S)(\uuu)\,d Q^{H^C_T\circ(\HHH^C_T)^\leftarrow}(\uuu)                           \\
&  = &  \int_{\I^d}          (H^C_T \circ S \circ (\HHH^C_T)^\leftarrow)(\uuu)\,d Q^{H^C_T\circ(\HHH^C_T)^\leftarrow}(\uuu)                           \\
&  = &  \int_{\overline\R^d} (H^C_T \circ S)(\xxx)                            \,d(Q^{H^C_T\circ(\HHH^C_T)^\leftarrow})_{(\HHH^C_T)^\leftarrow}(\xxx)  \\*
&  = &  \int_{\overline\R^d} (H^C_T \circ S)(\xxx)                            \,d Q^{H^C_T}                                                   (\xxx)  
\end{eqnarray*}
as was to be shown. 
\end{proof}

\bigskip
The previous results are remarkable since   $ C_T  $   and the restriction of   $ H^C_T $   to   $ \I^d $   are usually distinct.


\section{The Order Transform and Kendall's Tau}
\label{ot-kendall}

Throughout this section, 
we consider a fixed copula   $ C $. 
The discussion of Kendall's tau for the order transform   $ C_T $   of   $ C $   basically relies on the discussion of   $ [C_T,C_T] $. 

\bigskip
The following result provides a comparison of the copulas   $ C $   and   $ C_T $   in terms of Kendall's distribution function: 

\bigskip
\begin{theorem}{}
\label{t.order}
$ C_T $   satisfies 
$ K_{C_T} \leq K_C $. 
In particular, 
$$  [C,C]  \leq  [C_T,C_T]  $$%
\end{theorem}%

\bigskip
\begin{proof}
Lemma 
\ref{l-image} yields 
$$  (Q^{C_T})_{C_T}            
 =  ( Q^{H^C_T} )_{H^C_T}      
 =  ((Q^{H^C})_T)_{H^C_T}      
 =  (Q^{H^C})_{H^C_T \circ T}  
$$
and for every   $ \uuu\in\I^d $   we have 
$$    C(\uuu)
  =   Q^{H^C} \Bigl[ \Bigl\{ \vvv\in\I^d \Bigm|   \vvv  \leq   \uuu  \Bigr\} \Bigr]
\leq  Q^{H^C} \Bigl[ \Bigl\{ \vvv\in\I^d \Bigm| T(\vvv) \leq T(\uuu) \Bigr\} \Bigr]
  =   (H^C_T \circ T)(\uuu)
$$
For every   $ t\in\I^d $, 
we thus obtain 
\begin{eqnarray*}
        K_{C_T}(t)                                                                         
&  = &  (Q^{C_T})_{C_T}[[0,t]]    	                                                       \\[1ex]
&  = &  (Q^{H^C})_{H^C_T \circ T}[[0,t]]                                                   \\[.5ex]
&  = &  Q^C \Bigl[ \Bigl\{ \uuu\in\I^d \Bigm| (H^C_T \circ T)(\uuu) \leq t \Bigr\} \Bigr]  \\
&\leq&  Q^C \Bigl[ \Bigl\{ \uuu\in\I^d \Bigm|              C (\uuu) \leq t \Bigr\} \Bigr]  \\*[1ex]
&  = &  K_C(t)                                                                             
\end{eqnarray*}
which proves the first inequality. 
The second inequality then follows from Lemma 
\ref{distribution.l-Kendall-2}. 
\end{proof}

\bigskip
The next result provides several useful representations of   $ [C_T,C_T] $: 

\bigskip
\begin{theorem}{} 
\label{t.biconvex-1}
$ C_T $   satisfies 
\begin{eqnarray*}
        [C_T,C_T]                                                                           
&  = &  \int_{\I^d} Q^C \Bigl[                                                              
        \Bigl\{ \vvv\in\I^d \Bigm| T(\vvv) \leq        T(\uuu) \Bigr\} \Bigr] \,dQ^C(\uuu)  \\
&  = &  \int_{\I^d} Q^C\biggl[ \bigcup\nolimits_{\tilde\pi\in\tilde\Gamma^\pi}              
        \Bigl\{ \vvv\in\I^d \Bigm|   \vvv  \leq\tilde\pi(\uuu) \Bigr\} \biggr]\,dQ^C(\uuu)  \\*
&  = &  \int_{\I^d} H^C_T(T(\uuu))                                            \,dQ^C(\uuu)  
\end{eqnarray*}%
\end{theorem}%

\bigskip
\begin{proof}
Lemma 
\ref{l-image} yields 
\begin{eqnarray*}
        [C_T,C_T]                                                                                                                        
&  = &  \int_{\I^d}               C_T(\uuu)                                                                          \,dQ^{  C_T}(\uuu)  \\
&  = &  \int_{\overline\R^d}    H^C_T(\xxx)                                                                          \,dQ^{H^C_T}(\xxx)  \\
&  = &  \int_{\overline\R^d} Q^{H^C_T}\Bigl[ \Bigl\{ \sss\in\overline\R^d \Bigm|   \sss  \leq   \xxx  \Bigr\} \Bigr] \,dQ^{H^C_T}(\xxx)  \\
&  = &  \int_{\overline\R^d} Q^{H^C}  \Bigl[ \Bigl\{ \yyy\in\overline\R^d \Bigm| T(\yyy) \leq   \xxx  \Bigr\} \Bigr] \,dQ^{H^C_T}(\xxx)  \\
&  = &  \int_{\overline\R^d} Q^{H^C}  \Bigl[ \Bigl\{ \yyy\in\overline\R^d \Bigm| T(\yyy) \leq T(\zzz) \Bigr\} \Bigr] \,dQ^{H^C}  (\zzz)  \\*
&  = &  \int_{\I^d}          Q^C      \Bigl[ \Bigl\{ \vvv\in         \I^d \Bigm| T(\vvv) \leq T(\uuu) \Bigr\} \Bigr] \,dQ^C      (\uuu)  
\end{eqnarray*}
Moreover, 
for any   $ \uuu,\vvv\in\I^d $   there exist some   $ \tilde\pi_\uuu, \tilde\pi_\vvv \in \tilde\Gamma^\pi $   such that   
$ T(\uuu) = \tilde\pi_\uuu(\uuu) $   and 
$ T(\vvv) = \tilde\pi_\vvv(\vvv) $.  
Since   $ \tilde\Gamma^\pi $   is a group, 
this yields   $ T(\vvv) \leq T(\uuu) $   if and only if there exists some   $ \tilde\pi\in\tilde\Gamma^\pi $   such that   $ \vvv\leq\tilde\pi(\uuu) $. 
This yields 
\begin{eqnarray*}
        [C_T,C_T]                                                                            
&  = &  \int_{\I^d} Q^C\Bigl[                                                                
        \Bigl\{ \vvv\in\I^d \Bigm| T(\vvv) \leq         T(\uuu) \Bigr\} \Bigr] \,dQ^C(\uuu)  \\*
&  = &  \int_{\I^d} Q^C\biggl[ \bigcup\nolimits_{\tilde\pi\in\tilde\Gamma^\pi}               
        \Bigl\{ \vvv\in\I^d \Bigm|   \vvv  \leq \tilde\pi(\uuu) \Bigr\} \biggr]\,dQ^C(\uuu)  
\end{eqnarray*}
Finally, 
we have 
\begin{eqnarray*}
        [C_T,C_T]                                                                                               
&  = &  \int_{\I^d} Q^C      \Bigl[ \Bigl\{ \vvv\in\I^d \Bigm| T(\vvv) \leq T(\uuu) \Bigr\} \Bigr] \,dQ^C(\uuu)  \\
&  = &  \int_{\I^d} Q^{H^C_T}\Bigl[ \Bigl\{ \www\in\I^d \Bigm|   \www  \leq T(\uuu) \Bigr\} \Bigr] \,dQ^C(\uuu)  \\*
&  = &  \int_{\I^d}    H^C_T (T(\uuu))                                                             \,dQ^C(\uuu)  
\end{eqnarray*}
which completes the proof. 
\end{proof}

\bigskip
\begin{corollary}{}
\label{c.biconvex}
\begin{onelist}
\item   If   $ [C_T,C_T] = 0 $, 
then $ [C,C] = 0 $.
\item   If   $ C $   is symmetric with   $ [C,C] = 0 $, 
then   $ [C_T,C_T] = 0 $. 
\end{onelist}
\end{corollary}

\bigskip
\begin{proof}
Because of Theorem 
\ref{t.order} we have   $ [C,C] \leq [C_T,C_T] $,   
which yields (1). 
From Theorem 
\ref{t.biconvex-1} we obtain 
\begin{eqnarray*}
        [C_T,C_T]                                                                                                                                              
&  = &  \int_{\I^d} Q^C\biggl[ \bigcup\nolimits_{\tilde\pi\in\tilde\Gamma^\pi} \Bigl\{ \vvv\in\I^d \Bigm| \vvv\leq\tilde\pi(\uuu) \Bigr\} \biggr]\,dQ^C(\uuu)  \\
&\leq&  \sum_{\tilde\pi\in\tilde\Gamma^\pi} \int_{\I^d} Q^C\Bigl[              \Bigl\{ \vvv\in\I^d \Bigm| \vvv\leq\tilde\pi(\uuu) \Bigr\} \Bigr ]\,dQ^C(\uuu)  \\
&  = &  \sum_{\tilde\pi\in\tilde\Gamma^\pi} \int_{\I^d}  (C \circ \tilde\pi)(\uuu)                                                               \,dQ^C(\uuu)  \\
&  = &  \sum_{      \pi\in      \Gamma^\pi} \int_{\I^d}  (\pi(C))(\uuu)                                                                          \,dQ^C(\uuu)  \\*
&  = &  \sum_{      \pi\in      \Gamma^\pi} [\pi(C),C]                                                                                                         
\end{eqnarray*}
If   $ C $   is symmetric, 
then the previous inequality becomes   $ [C_T,C_T] \leq d!\,[C,C] $   since   $ |\Gamma^\pi| = d! $\,. 
This yields (2). 
\end{proof}

\bigskip
The following result resumes Theorem 
\ref{t.order} and Corollary 
\ref{c.biconvex} in terms of Kendall's tau: 

\bigskip
\begin{theorem}{}
\label{t.Kendall}
Kendall's tau satisfies 
$$  \kappa[C]  \leq  \kappa[C_T]  $$
In particular: 
\begin{onelist}
\item   If   $ C_T $   minimizes Kendall's tau, 
then   $ C $   minimizes Kendall's tau as well. 
\item   If   $ C $   is symmetric and minimizes Kendall's tau, 
then   $ C_T $   minimizes Kendall's tau as well. 
\end{onelist}
\end{theorem}
\pagebreak

\bigskip
We note in passing that Theorem 
\ref{t.Kendall} yields two results due to 
Dietz et al.\ [2016; Examples 4.2 and 4.3]: 

\bigskip
\begin{corollary}{}
\label{c.bounds}
The upper Fr{\'e}chet--Hoeffding bound satisfies   $ M_T=M $, 
and in the bivariate case 
the lower Fr{\'e}chet--Hoeffding bound satisfies   $ W_T=W $. 
\end{corollary}

\bigskip
\begin{proof}
It has been shown by 
Fuchs et al.\ [2018; Theorems 3.2 and 3.3] 
that the upper Fr{\'e}chet--Hoeffding bound is the only copula maximizing Kendall's tau and 
that in the bivariate case the lower Fr{\'e}chet--Hoeffding bound is the only copula minimizing Kendall's tau. 
The assertion follows. 
\end{proof}

\bigskip
Theorem 
\ref{t.Kendall} provides conditions under which    $ \kappa[C] = \kappa[C_T] $. 
Below we shall show that the product copula satisfies   $ \kappa[\Pi] < \kappa[\Pi_T] $; 
see Theorem 
\ref{t.product}. 

\bigskip
The following example provides a copula   $ D $   for which   $ D \not\leq_c D_T $; 
it thus shows that Theorem 
\ref{t.Kendall} cannot be obtained from the fact that Kendall's tau is monotone with respect to the concordance order   $ \leq_c $: 

\bigskip
\begin{example}{}
\label{e.bounds}
Assume that   $ d=2 $   and consider the copula 
\begin{eqnarray*}
        D              
& := &  \frac{M+W}{2}  
\end{eqnarray*}
Let   $ F_{(0,1)} $   denote the distribution function corresponding to the uniform distribution on   $ (0,1) $. 
According to 
Dietz et al.\ [2016; Examples 4.2 and 4.3], 
we have 
\begin{eqnarray*}
        H^D_T(\xxx)                                                                              
&  = &  \frac{1}{2}\,H^M_T(\xxx) + \frac{1}{2}\,H^W_T(\xxx)                                      \\*
&  = &  \frac{1}{2} \, F_{(0,1)}(\min\{x_1,x_2\})                                                
      + \frac{1}{2} \, \Bigl( F_{(0,1)}(2x_1) + \Bigl( 2\,F_{(0,1)}(x_2)-1 \Bigr)^+ -1 \Bigr)^+  
\end{eqnarray*}
and hence 
\begin{eqnarray*}
        \HHH^D_T(\xxx)                                                                         
&  = &  \leftmatrix{c}                                                                         
		 \dps \frac{1}{2}\,F_{(0,1)} (x_1) + \frac{1}{2}\,F_{(0,1)}(2x_1)                      \\[2ex]
		 \dps \frac{1}{2}\,F_{(0,1)} (x_2) + \frac{1}{2}\,\Bigl( 2\,F_{(0,1)}(x_2)-1 \Bigr)^+  
		\rightmatrix                                                                           
\end{eqnarray*}
Since the univariate marginal distribution functions of   $ H^D_T $   are continuous and strictly increasing on   $ (0,1) $, 
we obtain 
$$  D                             
    \leftmatrix{c}                
       3/8                        \\
       4/8                        
      \rightmatrix                
 =  \frac{3}{16}                  
 >  \frac{2}{16}                  
 =  (H^D_T\circ(\HHH^D_T )^{-1})  
    \leftmatrix{c}                
     3/8                          \\
     4/8                          
    \rightmatrix                  
 =  D_T                           
    \leftmatrix{c}                
     3/8                          \\
     4/8                          
    \rightmatrix   
$$
and hence   $ D \not\leq D_T $. 
Since   $ d=2 $, 
the concordance order agrees with the pointwise order and we obtain   $ D \not\leq_c D_T $. 
\end{example}

\bigskip
The following example shows that the map   $ C \mapsto C_T $   is not order preserving with respect to the concordance order: 

\bigskip
\begin{example}{}
\label{e.order}
Assume that   $ d=2 $. 
For every symmetric copula   $ C : \I^2\to\I $, 
Theorem 
\ref{t.biconvex-1} yields 
\begin{eqnarray*}
        [C_T,C_T]                                                                 
&  = &  2\,[C,C] - \int_{\I^2} C(\min\{u_1,u_2\},\min\{u_1,u_2\})\,dQ^C(u_1,u_2)  
\end{eqnarray*}
and for the computation of   $ [C,C] $   in the case where   $ C $   is a shuffle we refer to 
Fuchs et al.\ [2018; Theorem 5.1]. 
The copula   $ A: \I^2\to\I $   
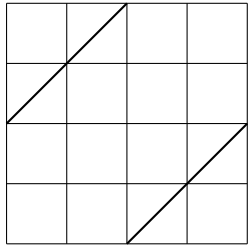
\begin{center}
\begin{tikzpicture}[xscale=0.8,yscale=0.8] 
\draw [thin]  (0,0) -- (0,4);
\draw [thin]  (1,0) -- (1,4);
\draw [thin]  (2,0) -- (2,4);
\draw [thin]  (3,0) -- (3,4);
\draw [thin]  (4,0) -- (4,4);
\draw [thin]  (0,0) -- (4,0);
\draw [thin]  (0,1) -- (4,1);
\draw [thin]  (0,2) -- (4,2);
\draw [thin]  (0,3) -- (4,3);
\draw [thin]  (0,4) -- (4,4);
\draw [thick] (0,2) -- (2,4);
\draw [thick] (2,0) -- (4,2);
\end{tikzpicture}
\end{center}
defined as the shuffle of   $ M $   with respect to the shuffling structure   $ \{[\aaa_i,\bbb_i]\}_{i\in\{1,\dots,4\}} $   with 
$$\begin{array}{ccccccc}
 \aaa_1 &=& (\phantom{00}0,2/4)  &&  \bbb_1 &=& (1/4,          3/4)  \\
 \aaa_2 &=& (1/4          ,3/4)  &&  \bbb_2 &=& (2/4,\phantom{00}1)  \\
 \aaa_3 &=& (2/4,\phantom{00}0)  &&  \bbb_3 &=& (3/4          ,1/4)  \\
 \aaa_4 &=& (3/4,          1/4)  &&  \bbb_4 &=& (\phantom{00}1,2/4)  
\end{array}$$
is symmetric and satisfies   $ [A_T,A_T] = 1/2 $,   
and the copula   $ B : \I^2\to\I $   
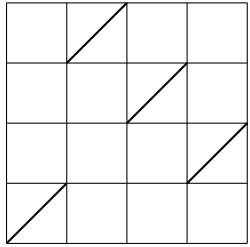
\begin{center}
\begin{tikzpicture}[xscale=0.8,yscale=0.8] 
\draw [thin]  (0,0) -- (0,4);
\draw [thin]  (1,0) -- (1,4);
\draw [thin]  (2,0) -- (2,4);
\draw [thin]  (3,0) -- (3,4);
\draw [thin]  (4,0) -- (4,4);
\draw [thin]  (0,0) -- (4,0);
\draw [thin]  (0,1) -- (4,1);
\draw [thin]  (0,2) -- (4,2);
\draw [thin]  (0,3) -- (4,3);
\draw [thin]  (0,4) -- (4,4);
\draw [thick] (0,0) -- (1,1);
\draw [thick] (1,3) -- (2,4);
\draw [thick] (2,2) -- (3,3);
\draw [thick] (3,1) -- (4,2);
\end{tikzpicture}
\end{center}
defined as the shuffle of   $ M $   with 
$$\begin{array}{ccccccc}
 \aaa_1 &=& (\phantom{00}0,\phantom{00}0)  &&  \bbb_1 &=& (1/4,          1/4)  \\
 \aaa_2 &=& (1/4          ,          3/4)  &&  \bbb_2 &=& (2/4,\phantom{00}1)  \\
 \aaa_3 &=& (2/4          ,          2/4)  &&  \bbb_3 &=& (3/4          ,3/4)  \\
 \aaa_4 &=& (3/4          ,          1/4)  &&  \bbb_4 &=& (\phantom{00}1,2/4)  
\end{array}$$
is symmetric as well and satisfies   $ [B_T,B_T] = 3/8 $. 
\\
It it obvious that   $ A \leq B $. 
On the other hand, 
we have   $ \kappa[A_T] = 1 $, 
and it then follows from 
Fuchs et al.\ [2018; Theorem 3.2] 
that   $ A_T = M $.
This yields   $ B_T \leq M = A_T $, 
and from   $ [B_T,B_T] < [A_T,A_T] $   we obtain   $ A_T \neq B_T $   and hence   $ A_T \not\leq B_T $. 
Since   $ d=2 $, 
this yields   $ A \leq_c B $   and   $ A_T \not\leq_c B_T $. 
\end{example}

\bigskip
Returning to the case of arbitrary dimension   $ d $, 
we conclude this section with another representation of   $ [C_T,C_T] $   for a class of copulas 
which includes every copula whose copula measure is absolutely continuous with respect to Lebesgue measure. 
This class includes the product copula, 
which will be studied in the following section. 

\bigskip
For   $ \tilde\pi\in\tilde\Gamma^\pi $, 
define
\begin{eqnarray*}
        A_{\tilde\pi}                                                                                                                       
& := &  \Bigl\{ \uuu \in \I^d \Bigm|  (\tilde\pi(\uuu))_i < (\tilde\pi(\uuu))_{i+1} \;\text{\rm for all}\; i\in\{1,\dots,d\!-\!1\} \Bigr\}  
\end{eqnarray*}
Then every    $ \uuu \in A_{\tilde\pi} $   satisfies   $ T(\uuu) = \tilde\pi(\uuu) $    
and the family   $ \fa{A}{{\tilde\pi}}{{\tilde\Gamma^\pi}} $   is disjoint. 

\bigskip
\begin{lemma}{} 
\label{l-biconvex-2}
Assume that   $ C $   satisfies   $ Q^C[ \sum_{\tilde\pi\in\tilde\Gamma^\pi} A_{\tilde\pi} ] = 1 $.
Then 
$$  [C_T,C_T]  =  \sum_{\pi\in\Gamma^\pi} \int_{T(\I^d)} H^C_T (\uuu) \,dQ^{\pi(C)}(\uuu)  $$
In particular, 
if   $ C $   is also symmetric, 
then 
$$  [C_T,C_T]  = d! \int_{T(\I^d)} H^C_T(\uuu) \,dQ^C(\uuu)  $$%
\end{lemma}%

\bigskip
\begin{proof}{}
Since   $ Q^C[ \sum_{\tilde\pi\in\tilde\Gamma^\pi} A_{\tilde\pi} ] = 1 $, 
Theorem 
\ref{t.biconvex-1} yields 
\begin{eqnarray*} 
        [C_T,C_T]                                                                                                                 
&  = &  \int_{\I^d}                                                                   H^C_T(        T(\uuu)) \,dQ^C       (\uuu)  \\
&  = &  \int_{\I^d} \sum_{\tilde\pi\in\tilde\Gamma^\pi} \chi_{A_{\tilde\pi}}(\uuu) \, H^C_T(        T(\uuu)) \,dQ^C       (\uuu)  \\
&  = &  \sum_{\tilde\pi\in\tilde\Gamma^\pi} \int_{\I^d} \chi_{A_{\tilde\pi}}(\uuu) \, H^C_T(\tilde\pi(\uuu)) \,dQ^C       (\uuu)  \\*
&  = &  \sum_{\pi\in\Gamma^\pi} \int_{T(\I^d)}                                        H^C_T(          \vvv ) \,dQ^{\pi(C)}(\vvv)  
\end{eqnarray*}
This proves the assertion.
\end{proof}


\section{The Order Transform of the Product Copula}
\label{ot-product}

In the present section we determine Kendall's tau for the order transform   $ \Pi_T $   of the product copula   $ \Pi $   
and for certain marginals of the order transform. 
To this end, 
we first recall a formula for the distribution function   $ H^\Pi_T $; 
see Dietz et al.\ [2016; Example 5.3]: 

\bigskip
\begin{lemma}{}
\label{l.product-identity}
The product copula   $ \Pi $   satisfies 
$$  H^\Pi_T(\uuu)  =  d!\,\det\Bigl[ (a_{i,j}(\uuu))_{i,j\in\{1,\dots,d\}} \Bigr]  $$
for every   $ \uuu \in T(\I^d) $, 
where 
$$  a_{i,j}(\uuu)                                                               
:=  \begin{cases}                                                               
     \dps \frac{u_i^{j-i+1}}{(j\!-\!i\!+\!1)!}  &   \text{if   $ i \leq j+1 $}  \\[2ex]
     0                                          &   \text{else}                 
    \end{cases}                                                                 
$$
for all   $ i,j\in\{1,\dots,d\} $. 
\end{lemma}%

\bigskip
We also recall that the copula measure of the product copula is the Lebesgue measure, 
which means that Lemma 
\ref{l-biconvex-2} applies to the product copula. 

\bigskip
\begin{theorem}{}
\label{t.product}
The product copula satisfies  
$$  [\Pi,\Pi]      
 =  \frac{1}{2^d}  
 <  \frac{1}{d+1}  
 =  [\Pi_T,\Pi_T]  
$$
and hence 
$$  \kappa[\Pi]                               
 =  0                                         
 <  \frac{2^d - (d\!+\!1)}{(2^{d-1}-1)(d+1)}  
 =  \kappa[\Pi_T]                             
$$
In particular, 
$ \lim_{d\to\infty} \kappa[\Pi_T] = 0 $. 
\end{theorem}

\bigskip
\begin{proof}{}
Since the product copula is absolutely continuous and symmetric, 
Lemma 
\ref{l-biconvex-2} together with Lemma 
\ref{l.product-identity} and Lemma 
\ref{l.app-1} yields 
\begin{eqnarray*} 
        [\Pi_T,\Pi_T]                                                                                    
&  = &  d! \int_{T(\I^d)} H^\Pi_T(\uuu)                                                 \,dQ^\Pi (\uuu)  \\
&  = &  d! \int_{T(\I^d)} d! \, \det\Bigl[ (a_{i,j}(\uuu))_{i,j\in\{1,\dots,d\}} \Bigr] \,d\leb^d(\uuu)  \\
&  = &  d!\,d!\,\frac{1}{d!\,(d+1)!}                                                                     \\*
&  = &  \frac{1}{d+1}                                                                                    
\end{eqnarray*}
The identity for   $ [\Pi,\Pi] $   is well known. 
\end{proof}

\bigskip
Similar but more complicated identities can be obtained for certain margins of the order transform of the product copula. 
In the sequel, 
we use the canonical extensions 
of transformations of   $ \I^d $   
to transformations of   $ \overline\R^d $   without changing the notation. 
We also need the following definitions: 

\bigskip
Consider   $ K\subseteq\{1,\dots,d\} $   such that   $ |K|\geq2 $.   
Then there exists a unique strictly increasing sequence   $ \fa{k}{j}{\{1,\dots,|K|\}} \subseteq \{1,\dots,d\} $   
such that   $ K = \fa{k}{j}{\{1,\dots,|K|\}} $   and we denote 
\begin{hylist}
\item   by   $ \overline\R^K $   a copy of   $ \overline\R^{|K|} $   with coordinates   $ k_1,\dots,k_{|K|} $   instead of   $ 1,\dots,|K| $, 
\item   by   $ \I^K $   the unit cube of   $ \overline\R^K $,   and 
\item   by   $ \CC^K $   the collection of all copulas   $ \I^K\to\I $. 
\end{hylist}
We also define a map   $ \tilde\varrho_K : \I^K\to\I^d $   by letting 
\begin{eqnarray*}
        \tilde\varrho_K(\vvv)                                                                         
& := &  \sum_{j \in \{1,\dots,|K|\}} v_j\,\eee_{k_j} + \sum_{k \in \{1,\dots,d\} \setminus K} \eee_k  
\end{eqnarray*}
and the map   $ \varrho_K : \CC^d\to\CC^K $   given by 
\begin{eqnarray*}
        \varrho_K(C)             
& := &  C \circ \tilde\varrho_K  
\end{eqnarray*}
Then   $ \varrho_K $   associates indeed with every copula in   $ \CC^d $   a copula in   $ \CC^K $, 
which is called the 
\emph{margin} of   $ C $   with respect to   $ K $. 
We have the following result: 

\bigskip
\begin{theorem}{}
\label{t.product-K}
$ \varrho_K(\Pi_T) $   satisfies 
\begin{eqnarray*}
        \Bigl[ \varrho_K(\Pi_T),\varrho_K(\Pi_T) \Bigr]                            
&  = &     \int_{\I^d}            \Pi_T(\eeta_K(\eins,\uuu)) \,dQ^{  \Pi_T}(\uuu)  \\*
&  = &     \int_{\overline\R^d} H^\Pi_T(\eeta_K(\eins,\xxx)) \,dQ^{H^\Pi_T}(\xxx)  \\*
&  = &  d! \int_{T(\I^d)}       H^\Pi_T(\eeta_K(\eins,\uuu)) \,dQ^{  \Pi  }(\uuu)  
\end{eqnarray*}%
\end{theorem}%

\bigskip
\begin{proof}
Consider the map   $ \eeta_{K,\eins} : \I^d\to\I^d $   given by   $ \eeta_{K,\eins}(\uuu) := \eeta_K(\eins,\uuu) $. 
Then we obtain 
\begin{eqnarray*}
        \Bigl[ \varrho_K(\Pi_T),\varrho_K(\Pi_T) \Bigr]                                                                                               
&  = &  \Bigl[ \Pi_T\circ\tilde\varrho_K,\Pi_T\circ\tilde\varrho_K \Bigr]                                                                             \\[.5ex]
&  = &  \int_{\I^K}                                        \Pi_T(\tilde\varrho_K(\vvv)) \,d Q^{\Pi_T \circ \tilde\varrho_K}                   (\vvv)  \\
&  = &  \int_{\tilde\varrho_K^{-1}(\eeta_{K,\eins}(\I^d))} \Pi_T(\tilde\varrho_K(\vvv)) \,d Q^{\Pi_T \circ \tilde\varrho_K}                   (\vvv)  \\
&  = &  \int_{\eeta_{K,\eins}(\I^d)}                       \Pi_T(\uuu)                  \,d(Q^{\Pi_T \circ \tilde\varrho_K})_{\tilde\varrho_K}(\uuu)  \\
&  = &  \int_{\eeta_{K,\eins}(\I^d)}                       \Pi_T(\uuu)                  \,d Q^{\Pi_T}                                         (\uuu)  \\
&  = &  \int_{\eeta_{K,\eins}(\I^d)}                       \Pi_T(\uuu)                  \,d(Q^{\Pi_T})_{\eeta_{K,\eins}}                      (\uuu)  \\
&  = &  \int_{\I^d}                                        \Pi_T(\eeta_{K,\eins}(\uuu)) \,d Q^{\Pi_T}                                         (\uuu)  \\*
&  = &  \int_{\I^d}                                        \Pi_T(\eeta_K (\eins, \uuu)) \,d Q^{\Pi_T}                                         (\uuu)  
\end{eqnarray*}
which gives the first identity. 
The second identity then follows from the first identity and Lemma 
\ref{l.commute}. 
Furthermore, 
we have 
\begin{eqnarray*}
        \Bigl[ \varrho_K(\Pi_T),\varrho_K(\Pi_T) \Bigr]                                                                        
&  = &  \int_{\overline\R^d}                                     H^\Pi_T(\eeta_K(\eins,  \xxx ))         \,dQ^{H^\Pi_T}(\xxx)  \\
&  = &  \int_{\overline\R^d}                                     H^\Pi_T(\eeta_K(\eins,T(\xxx)))         \,dQ^{H^\Pi}  (\xxx)  \\
&  = &  \int_{\I^d}                                              H^\Pi_T(\eeta_K(\eins,T(\uuu)))                               
        \Biggl( \sum_{\tilde\pi\in\tilde\Gamma^\pi} \chi_{A_{\tilde\pi}}(\uuu) \Biggr)                   \,dQ^\Pi      (\uuu)  \\
&  = &  \sum_{\tilde\pi\in\tilde\Gamma^\pi} \int_{A_{\tilde\pi}} H^\Pi_T(\eeta_K(\eins,\tilde\pi(\uuu))) \,dQ^\Pi      (\uuu)  \\*
&  = &  d! \int_{T(\I^d)}                                        H^\Pi_T(\eeta_K(\eins,          \uuu )) \,dQ^\Pi      (\uuu)  
\end{eqnarray*}
which gives the last identity. 
\end{proof}

\bigskip
From the previous result, 
combined with Lemma 
\ref{l.product-identity}, 
we obtain explicit formulas for   $ [\varrho_K(\Pi_T),\varrho_K(\Pi_T)] $,   
and hence for   $ \kappa[\varrho_K(\Pi_T)] $,   
in the case   $ K=\{1,\dots,k\} $: 

\bigskip
\begin{corollary}{}
\label{c.product-lowertail}
$ \varrho_{\{1,\dots,k\}}(\Pi_T) $   satisfies 
\begin{eqnarray*}
        \Bigl[ \varrho_{\{1,\dots,k\}}(\Pi_T),\varrho_{\{1,\dots,k\}}(\Pi_T) \Bigr]                                        
&  = &  \frac{1}{2} - \frac{1}{4} \sum_{h=2}^k \frac{1}{2h-1}\, \binom{2h}{h} \binom{2d+2-2h}{d+1-h} \bigg/ \binom{2d}{d}  \\*
&  = &  \frac{1}{d+1} + \frac{1}{2d} \sum_{l=1}^{d-k} \binom{d}{l-1} \binom{d}{l} \bigg/ \binom{2d-1}{2l-1}                
\end{eqnarray*}
for every   $ k\in\{2,\dots,d\} $. 
In particular, 
\begin{eqnarray*}
        \kappa[\varrho_{\{1,\dots,k\}}(\Pi_T)]                                                                                 
&  = &  1 - \frac{2^{k-2}}{2^{k-1}-1} \sum_{h=2}^k \frac{1}{2h-1}\, \binom{2h}{h} \binom{2d+2-2h}{d+1-h} \bigg/ \binom{2d}{d}  
\end{eqnarray*}
holds for every   $ k\in\{2,\dots,d\} $, 
the sequence   $ \{ \kappa[\varrho_{\{1,\dots,k\}}(\Pi_T)] \}_{k\in\{1,\dots,d\}} $   is decreasing with 
\begin{eqnarray*}
        \kappa[\varrho_{\{1,\dots,d\}}(\Pi_T)]  
&  = &  \frac{1}{d+1}                           
\end{eqnarray*} 
and for every   $ k\in\{2,3,\dots\} $   the sequence   $ \{ \kappa[\varrho_{\{1,\dots,k\}}(\Pi_T)] \}_{d\in\{k,k+1,\dots\}} $   is increasing with 
\begin{eqnarray*}
        \lim_{d\to\infty} \kappa[\varrho_{\{1,\dots,k\}}(\Pi_T)]                                        
&  = &  1 - \frac{2^{k-2}}{2^{k-1}-1} \sum_{h=2}^k \frac{1}{2h-1}\, \binom{2h}{h} \,\frac{1}{2^{2h-2}}  
\end{eqnarray*}%
\end{corollary}%

\bigskip
\begin{proof}
Since   $ \eeta_{\{1,\dots,k\}}(\eins,\uuu) \in T(\I^d) $   holds for every   $ \uuu \in T(\I^d) $, 
Theorem 
\ref{t.product-K} together with Lemma 
\ref{l.product-identity} and Corollary 
\ref{c.app-1} yields 
\begin{eqnarray*}
        \Bigl[ \varrho_{\{1,\dots,k\}}(\Pi_T),\varrho_{\{1,\dots,k\}}(\Pi_T) \Bigr]                                                                           
&  = &  d! \int_{T(\I^d)}                                        H^\Pi_T(\eeta_{\{1,\dots,k\}}(\eins,\uuu))                             \,dQ^\Pi (\uuu)       \\
&  = &  d! \int_{T(\I^d)} d!                       \,\det\Bigl[ (a_{i,j}(\eeta_{\{1,\dots,k\}}(\eins,\uuu)))_{i,j\in\{1,...,d\}} \Bigr] \,d\leb^d(\uuu)       \\
&  = &  (d!)^2 \biggl( \frac{1}{d!\,(d\!+\!1)!} + \frac{1}{(2d)!} \sum_{l=1}^{d-k} \frac{(2d\!-\!2l)!}{(d\!-\!l)!\,(d\!+\!1\!-\!l)!} \binom{2l-1}{l} \biggr)  \\*
&  = &  \frac{1}{d+1}                      + \frac{1}{2d} \sum_{l=1}^{d-k} \binom{d}{l-1} \binom{d}{l} \bigg/ \binom{2d-1}{2l-1}                              
\end{eqnarray*}
which gives the second identity. 
On the other hand, 
Lemma 
\ref{l.combi} yields 
\begin{eqnarray*} 
	    \sum_{l=1}^{d-k} \frac{(2d-2l)!}{(d-l)!\,(d+1-l)!}\,\frac{(2l-1)!}{l! (l-1)!}             
&  = &  \frac{1}{4} \sum_{l=1}^{d-k} \binom{2d+2-2l}{d+1-l} \binom{2l}{l} \frac{1}{2d-2l+1}       \\
&  = &  \frac{1}{4} \sum_{h=k+1}^d   \binom{2h}{h} \binom{2d+2-2h}{d+1-h} \frac{1}{2h-1}   	      \\
&  = &  \frac{1}{4} \biggl(                      - \binom{2d+2}{d+1}      \frac{1}{2d+1}          
                  - \sum_{h=0}^k     \binom{2h}{h} \binom{2d+2-2h}{d+1-h} \frac{1}{2h-1} \biggr)  \\*
&  = &  \frac{1}{4} \biggl( \binom{2d}{d} \frac{2(d-1)}{d+1}                                      
				  - \sum_{h=2}^k     \binom{2h}{h} \binom{2d+2-2h}{d+1-h} \frac{1}{2h-1} \biggr)  
\end{eqnarray*}
and hence 
\begin{eqnarray*}  
        \Bigl[ \varrho_{\{1,\dots,k\}}(\Pi_T),\varrho_{\{1,\dots,k\}}(\Pi_T) \Bigr]                              
&  = &  \frac{1}{d+1}                                                                                            
      + \frac{(d!)^2}{(2d)!} \sum_{l=1}^{d-k} \frac{(2d\!-\!2l)!}{(d\!-\!l)!\,(d\!+\!1\!-\!l)!} \binom{2l-1}{l}  \\
&  = &  \frac{1}{d+1}                                                                                            
      + \frac{1}{2}\,\frac{d-1}{d+1}                                                                             
      - \frac{1}{4}\,\frac{(d!)^2}{(2d)!} \sum_{h=2}^k \binom{2h}{h} \binom{2d+2-2h}{d+1-h} \frac{1}{2h-1}       \\*
&  = &  \frac{1}{2}                                                                                              
      - \frac{1}{4} \sum_{h=2}^k \frac{1}{2h-1}\, \binom{2h}{h} \binom{2d+2-2h}{d+1-h} \bigg/ \binom{2d}{d}      
\end{eqnarray*}
which gives the first identity. 
\end{proof}

\bigskip
The first identity of Corollary 
\ref{c.product-lowertail} is suitable for small   $ k $   while the second is suitable for large  $ k $. 
For   $ k=2 $, 
Corollary 
\ref{c.product-lowertail} yields 
\begin{eqnarray*}
        \kappa[\varrho_{\{1,2\}}(\Pi_T)]  
&  = &  \frac{d-1}{2d-1}                  
\end{eqnarray*}
which is in accordance with the literature; 
see 
Av{\'e}rous et al.\ [2005] and 
Navarro and Balakrishnan [2010]. 
In particular, 
the sequence   $ \{ \kappa[\varrho_{\{1,2\}}(\Pi_T)] \}_{d\in\{2,3,\dots\}} $   is increasing 
with   $ \lim_{d\to\infty}\kappa[\varrho_{\{1,2\}}(\Pi_T)] = 1/2 $. 
We illustrate Corollary 
\ref{c.product-lowertail} by numerical examples for certain values of   $ d $   and   $ k\in\{2,\dots,d\} $: 

\bigskip
\begin{examples}{}
\label{e.product-lowertail}
For   $ d\in\{2,\dots,5\} $   and   $ k\in\{2,\dots,d\} $, 
the table 
$$\begin{array}{|c|@{\;}cccc|}
\hline\rule{0ex}{2.5ex}%
   & \multicolumn{4}{c|}{d} \\
 k &    2    &    3    &    4    &    5    \\
\hline\rule{0ex}{2.5ex}%
 2 & 315/945 & 378/945 & 405/945 & 420/945 \\
 3 &         & 315/945 & 369/945 & 395/945 \\
 4 &         &         & 297/945 & 345/945 \\
 5 &         &         &         & 273/945 \\
\hline
\end{array}$$
presents the values of   $ \kappa[\varrho_{\{1,\dots,k\}}(\Pi_T)] $   for   $ d\in\{2,\dots,5\} $   and   $ k\in\{2,\dots,d\} $. 
\end{examples}

\bigskip
We now introduce a general reflection principle which, 
when combined with Corollary 
\ref{c.product-lowertail}, 
yields explicit formulas for   $ [\varrho_K(\Pi_T),\varrho_K(\Pi_T)] $, 
and hence for   $ \kappa[\varrho_K(\Pi_T)] $,   
in the case   $ K=\{d-k+1,\dots,d\} $. 
We shall need the following lemma: 

\bigskip
\begin{lemma}{}
\label{l.product-taurho}
The identity 
$$  \rho_K(\tau(\Pi_T))  =  \tau(\rho_K(\Pi_T))  $$
holds for every   $ K\subseteq\{1,\dots,d\} $   such that   $|K|\geq2 $. 
\end{lemma}

\bigskip
In the identity of Lemma 
\ref{l.product-taurho}, 
the transformation   $ \tau $   acts on   $ \CC^d $   on the left hand side and on   $ \CC^K $   on the right hand side. 
\pagebreak

\bigskip
\begin{proof}
According to 
Fuchs [2014; Theorem 4.1], 
every copula   $ D $   satisfies 
\begin{eqnarray*}
        (\tau(D))(\vvv)                                                                  
&  = &  \sum_{L \subseteq \{1,...,d\}} (-1)^{d-|L|} \, D(\eeta_L(\eins\!-\!\vvv,\eins))  
\end{eqnarray*}
For every   $ \uuu\in\I^K $   we obtain 
\begin{eqnarray*}
        (\rho_K(\tau(\Pi_T)))(\uuu)                                                                                                    
&  = &  (\tau(\Pi_T))(\tilde\rho_K(\uuu))                                                                                              \\[1ex]
&  = &  \sum_{L \subseteq \{1,...,d\}} (-1)^{d-|L|}  \, \Pi_T(\eeta_L                           (\eins\!-\!\tilde\rho_K(\uuu),\eins))  \\
&  = &  \sum_{L \subseteq \{1,...,d\}} (-1)^{d-|L|}  \, \Pi_T(\eeta_{\{1,\dots,d\} \setminus L} (\eins,\eins\!-\!\tilde\rho_K(\uuu)))  \\
&  = &  \sum_{J \subseteq \{1,...,d\}} (-1)^{|J|}    \, \Pi_T(\eeta_J                           (\eins,\eins\!-\!\tilde\rho_K(\uuu)))  \\
&  = &  \sum_{J \subseteq K}           (-1)^{|J|}    \, \Pi_T(\eeta_J                           (\eins,\eins\!-\!\tilde\rho_K(\uuu)))  \\
&  = &  \sum_{J \subseteq K}           (-1)^{|J|}    \, \Pi_T(\eeta_{\{1,\dots,d\} \setminus L} (\eins\!-\!\tilde\rho_K(\uuu),\eins))  \\
&  = &  \sum_{L \subseteq K}           (-1)^{|K|-|L|}\, \Pi_T(\eeta_L              (\tilde\rho_K(\eins\!-\!\uuu),\eins))               \\
&  = &  \sum_{L \subseteq K}           (-1)^{|K|-|L|}\, (\Pi_T \circ \tilde\rho_K) (\eeta_L     (\eins\!-\!\uuu ,\eins))               \\
&  = &  \sum_{L \subseteq K}           (-1)^{|K|-|L|}\, (\rho_K(\Pi_T))            (\eeta_L     (\eins\!-\!\uuu ,\eins))               \\*[1ex]
&  = &  (\tau(\rho_K(\Pi_T)))(\uuu)                                                                                                    
\end{eqnarray*}
Note that in the last three lines the transformation   $ \eeta_L $   acts on   $ \I^K $. 
\end{proof}

\bigskip
Consider now the map   $ b : \{1,\dots,d\}\to\{1,\dots,d\} $   given by 
\begin{eqnarray*}
        b(i)   
& := &  d-i+1  
\end{eqnarray*}
and the linear map   $ B : \overline\R^d\to\overline\R^d $   given by the $ (d \times d) $--matrix with entries 
\begin{eqnarray*}
        b_{i,j}                         
& := &  \begin{cases}                   
         1  &  \text{if   $ i+j=d+1 $}  \\
         0  &  \text{else}              
        \end{cases}                     
\end{eqnarray*}
We can now state the announced reflection principle: 

\bigskip
\begin{theorem}{}
\label{t.product-reflection}
The identity 
$$  \Bigl[ \varrho_{b(K)}(\Pi_T) , \varrho_{b(K)}(\Pi_T) \Bigr]  
 =  \Bigl[ \varrho_K(\Pi_T) , \varrho_K(\Pi_T) \Bigr]            
$$
holds for every   $ K\subseteq\{1,\dots,d\} $   such that   $ |K|\geq2 $. 
In particular, 
$$  \kappa[\varrho_{b(K)}(\Pi_T)]  =  \kappa[\varrho_K(\Pi_T)]  $$
holds for every   $ K\subseteq\{1,\dots,d\} $   such that   $ |K|\geq2 $. 
\end{theorem}
\pagebreak

\bigskip
\begin{proof}
Using Lemma 
\ref{l.product-taurho} and Theorem 
\ref{t.product-K} we obtain 
\begin{eqnarray*}
        \Bigl[      \varrho_{b(K)}(\Pi_T)  ,      \varrho_{b(K)}(\Pi_T)  \Bigr]       
&  = &  \Bigl[ \tau(\varrho_{b(K)}(\Pi_T)) , \tau(\varrho_{b(K)}(\Pi_T)) \Bigr]       \\
&  = &  \Bigl[ \varrho_{b(K)}(\tau(\Pi_T)) , \varrho_{b(K)}(\tau(\Pi_T)) \Bigr]       \\*
&  = &  \int_{\I^d} (\tau(\Pi_T))(\eeta_{b(K)}(\eins,\uuu)) \,dQ^{\tau(\Pi_T)}(\uuu)  
\end{eqnarray*}
and from Theorem 
\ref{t.product-K} we also obtain 
\begin{eqnarray*}
        \Bigl[ \varrho_K(\Pi_T) , \varrho_K(\Pi_T) \Bigr]                       
&  = &  \int_{\I^d}            \Pi_T(\eeta_K(\eins,\uuu)) \,dQ^{  \Pi_T}(\uuu)  \\*
&  = &  \int_{\overline\R^d} H^\Pi_T(\eeta_K(\eins,\xxx)) \,dQ^{H^\Pi_T}(\xxx)  
\end{eqnarray*}
It remains to show that the two integrals are identical. 
We have   $ B \circ T = \tilde\tau \circ T \circ \tilde\tau $   and hence   $ \tilde\tau \circ T = B \circ T \circ \tilde\tau $, 
and we also have   
$ Q^{H^\Pi_{T \circ \tilde\tau}} = (Q^{H^{\Pi}})_{T \circ \tilde\tau} = ((Q^{H^\Pi})_{\tilde\tau})_T = (Q^{H^\Pi})_T = Q^{H^\Pi_T} $. 
From these identities and Lemma 
\ref{l.commute} we obtain 
\begin{eqnarray*}
\lefteqn{                                                                                                              
 		\int_{\I^d}             \Pi_T (\eeta_L(\eeta_{b(K)}(\zero,  \uuu ),\eins))                                     
        \,dQ^{\Pi_T}  (\uuu)                                                                                          }\\*
&  = &  \int_{\overline\R^d}  H^\Pi_T (\eeta_L(\eeta_{b(K)}(\zero,  \xxx ),\eins))                                     
        \,dQ^{H^\Pi_T}(\xxx)                                                                                           \\
&  = &  \int_{\overline\R^d}  H^\Pi_T (\eeta_L(\eeta_{b(K)}(\zero,T(\xxx)),\eins))                                     
        \,dQ^{H^\Pi}  (\xxx)                                                                                           \\
&  = &  \int_{\overline\R^d}  H^\Pi_T (\eeta_L(\eins - \eeta_{b(K)}(\eins,(        \tilde\tau \circ T)(\xxx)),\eins))  
        \,dQ^{H^\Pi}  (\xxx)                                                                                           \\
&  = &  \int_{\overline\R^d}  H^\Pi_T (\eeta_L(\eins - \eeta_{b(K)}(\eins,(B \circ T \circ \tilde\tau)(\xxx)),\eins))  
        \,dQ^{H^\Pi}  (\xxx)                                                                                           \\
&  = &  \int_{\overline\R^d}  H^\Pi_T (\eeta_L(\eins - \eeta_{b(K)}(\eins,B(\zzz))                           ,\eins))  
        \,dQ^{H^\Pi_{T \circ \tilde\tau}}(\zzz)                                                                        \\
&  = &  \int_{\overline\R^d}  H^\Pi_T                            (\eeta_L(\eins - B \circ \eeta_K(\eins,\zzz),\eins))  
        \,dQ^{H^\Pi_T}                   (\zzz)                                                                        \\
&  = &  \int_{\overline\R^d}  H^\Pi_T                            ((B \circ \eeta_L)(\eins-\eeta_K(\eins,\zzz),\eins))  
        \,dQ^{H^\Pi_T}                   (\zzz)                                                                        \\
&  = &  \int_{\overline\R^d}  H^\Pi_{B \circ T}                           (\eeta_L (\eins-\eeta_K(\eins,\zzz),\eins))  
        \,dQ^{H^\Pi_T}                   (\zzz)                                                                        \\*
&  = &  \int_{\overline\R^d}  H^\Pi_{\tilde\tau \circ T \circ \tilde\tau} (\eeta_L (\eeta_K(\zero,\eins-\zzz),\eins))  
        \,dQ^{H^\Pi_T}                   (\zzz)                                                                        
\end{eqnarray*}
Using the first identity in the proof of Lemma 
\ref{l.product-taurho} we thus obtain 
\begin{eqnarray*}
\lefteqn{
        \Bigl[ \varrho_{b(K)}(\Pi_T),\varrho_{b(K)}(\Pi_T) \Bigr]                                                                                }\\*
&  = &  \int_{\I^d} (\tau(\Pi_T))(\eeta_{b(K)}(\eins,\uuu))                                                             \,dQ^{\tau(\Pi_T)}(\uuu)  \\
&  = &  \int_{\I^d} \sum_{L\subseteq\{1,\dots,d\}} (-1)^{d-|L|}\,                                                                                 
        \Pi_T            (\eeta_L(\eins-\eeta_{b(K)}(\eins,\uuu),\eins))                                                \,dQ^{\tau(\Pi_T)}(\uuu)  \\
&  = &  \sum_{L\subseteq\{1,\dots,d\}} (-1)^{d-|L|}\,\int_{\I^d}                                                                                  
        \Pi_T            (\eeta_L(\eeta_{b(K)}(\zero,\tilde\tau(\uuu)),\eins))                                          \,dQ^{\tau(\Pi_T)}(\uuu)  \\
&  = &  \sum_{L\subseteq\{1,\dots,d\}} (-1)^{d-|L|}\,\int_{\I^d}                                                                                  
        \Pi_T            (\eeta_L(\eeta_{b(K)}(\zero,\uuu),\eins))                                                      \,dQ^{\Pi_T}      (\uuu)  \\
&  = &  \sum_{L\subseteq\{1,\dots,d\}} (-1)^{d-|L|}\,\int_{\overline\R^d}                                                                         
        H^\Pi_{\tilde\tau \circ T \circ \tilde\tau} (\eeta_L(\eeta_K(\zero,\eins\!-\!\xxx),\eins))                      \,dQ^{H^\Pi_T}    (\xxx)  \\
&  = &  \int_{\overline\R^d} \sum_{L\subseteq\{1,\dots,d\}} (-1)^{d-|L|}\,                                                                        
                    Q^\Pi_{\tilde\tau \circ T \circ \tilde\tau}                                                                                   
        \Bigl[ \Bigl[ \zero,                         \eeta_L(\eeta_K(\zero,\eins\!-\!\xxx),\eins)         \Bigr] \Bigr] \,dQ^{H^\Pi_T}    (\xxx)  \\
&  = &  \int_{\overline\R^d} Q^\Pi_{\tilde\tau \circ T \circ \tilde\tau}                                                                          
        \Bigl[ \Bigl[                                        \eeta_K(\zero,\eins\!-\!\xxx)        , \eins \Bigr] \Bigr] \,dQ^{H^\Pi_T}    (\xxx)  \\
&  = &  \int_{\overline\R^d} Q^\Pi_{                 T \circ \tilde\tau}                                                                          
        \Bigl[ \Bigl[                        \zero , \eins - \eeta_K(\zero,\eins\!-\!\xxx)                \Bigr] \Bigr] \,dQ^{H^\Pi_T}    (\xxx)  \\
&  = &  \int_{\overline\R^d} Q^\Pi_T                                                                                                              
        \Bigl[ \Bigl[                        \zero ,         \eeta_K(\eins,          \xxx)                \Bigr] \Bigr] \,dQ^{H^\Pi_T}    (\xxx)  \\
&  = &  \int_{\overline\R^d} H^\Pi_T(                        \eeta_K(\eins,          \xxx))                             \,dQ^{H^\Pi_T}    (\xxx)  \\*
&  = &  \Bigl[ \varrho_K(\Pi_T),\varrho_K(\Pi_T) \Bigr]                                                                                           
\end{eqnarray*}
as was to be shown. 
\end{proof}

\bigskip
Theorem 
\ref{t.product-reflection} extends a result of 
Av{\'e}rous et al.\ [2005; Proposition 10] for   $ |K|=2 $. 
As noted before, 
combining Theorem 
\ref{t.product-reflection} and Corollary  
\ref{c.product-lowertail} yields explicit formulas for   $ [\varrho_K(\Pi_T),\varrho_K(\Pi_T)] $,   
and hence for   $ \kappa[\varrho_K(\Pi_T)] $,   
in the case   $ K=\{d-k+1,\dots,d\} $. 

\bigskip
As a final remark, 
we note that Theorem 
\ref{t.product-K} combined with Lemma 
\ref{l.product-identity} provides a general tool to compute Kendall's tau for the margins of   $ \Pi_T $   with respect to any   $ K\subseteq\{1,\dots,d\} $. 
We illustrate this by the following example which is not covered by Corollary  
\ref{c.product-lowertail}: 

\bigskip
\begin{example}{}
Assume that   $ d=5 $   and consider   $ K=\{1,2,3,5\} $. 
Since 
\linebreak
  $ Q^\Pi_T[ \I^d \setminus T(\I^d) ] = 0 $, 
Theorem 
\ref{t.product-K} together with Lemma 
\ref{l.product-identity} yields 
\begin{eqnarray*}
        \Bigl[ \rho_{\{1,2,3,5\}}(\Pi_T), \rho_{\{1,2,3,5\}}(\Pi_T) \Bigr]  
&  = &  5! \int_{T(\I^5)} H^\Pi_T((u_1,u_2,u_3,  1,u_5)) \,dQ^\Pi(\uuu)     \\
&  = &  5! \int_{T(\I^5)} H^\Pi_T((u_1,u_2,u_3,u_5,u_5)) \,dQ^\Pi(\uuu)     \\
&  = &  (5!)^2 \int_{T(\I^5)} \det                                          
		\left( \begin{matrix}                                               
		  u_1 & u_1^2/2! & u_1^3/3! & u_1^4/4! & u_1^5/5!                   \\
		    1 &      u_2 & u_2^2/2! & u_2^3/3! & u_2^4/4!                   \\
			0 &        1 &      u_3 & u_3^2/2! & u_3^3/3!                   \\
			0 &        0 &        1 &      u_5 & u_5^2/2!                   \\
			0 &        0 &        0 &        1 &      u_5                   
		\end{matrix} \right)                               dQ^\Pi(\uuu)     \\
&  = &  (5!)^2 \frac{47}{3\,628\,800}                                       \\*
&  = &  \frac{47}{252}                                                      
\end{eqnarray*}
We thus obtain 
\begin{eqnarray*}
        \kappa[\rho_{\{1,2,3,5\}}(\Pi_T)]                                               
&  = &  \frac{2^4\,[\rho_{\{1,2,3,5\}}(\Pi_T),\rho_{\{1,2,3,5\}}(\Pi_T)]-1}{2^{4-1}-1}  
\,\;=\;\, \frac{125}{441}                                                               
\end{eqnarray*}
and Theorem 
\ref{t.product-reflection} then yields   $ \kappa[\rho_{\{1,3,4,5\}}(\Pi_T)] = \kappa[\rho_{\{1,2,3,5\}}(\Pi_T)] = 125/441 $. 
\end{example}


\appendix

\section{Auxiliary Results}
\label{aux-results}

For a family   $ \{a_{i,j}\}_{i,j\in\N} $   of real numbers satisfying 
$ a_{j+1,j}=1 $   for every   $ j\in\N $   and 
$ a_{j+k,j}=0 $   for every   $ j\in\N $   and all   $ k\geq2 $, 
we define a sequence of matrices   $ \fa{A}{n}{\N_0} $   by letting   
$ A_0 := (1) $   and   
$ A_n := (a_{i,j})_{i,j\in\{1,\dots,n\}} $   for every   $ n\in\N $. 
Then the identity 
\begin{eqnarray*}
        \det(A_n)                                      
&  = &  \sum_{i=1}^n (-1)^{i+n} a_{i,n} \det(A_{i-1})  
\end{eqnarray*}
holds for every   $ n\in\N $; 
see 
Cahill et al.\ [2002]. 

\bigskip
Let   $ S_0 := 1 $. 
For   $ n\in\N $   we define a map $ S_n: \I^n\to\R $   by letting 
\begin{eqnarray*}
        S_n (\uuu)                                                                      
& := &  \sum_{i=1}^n (-1)^{i+n} \frac{u_i^{n-i+1}}{(n-i+1)!}\,S_{i-1}(u_1,...,u_{i-1})  
\end{eqnarray*}
Then   $ S_n(\uuu) $   is the determinant of the matrix   $ A_n(\uuu) $   with entries 
\begin{eqnarray*}
        a_{i,j}(\uuu)                                                  
& := &  \begin{cases}                                                  
         \dps \frac{u_i^{j-i+1}}{(j-i+1)!}  & \text{if}\;  i \leq j+1  \\*[2ex]
         \dps 0                             & \text{else}              
        \end{cases}                                                    
\end{eqnarray*}
for   $ n\in\N $   and   $ i,j\in\{1,\dots,n\} $. 

\bigskip
For the convenience of the reader, 
we recall some combinatorial identities: 

\bigskip
\begin{lemma}{}
\label{l.combi}
The identities 
\begin{eqnarray*}
        \sum_{k=0}^n \frac{1}{2k-1} \binom{2n-2k}{n-k} \binom{2k}{k}       
&  = &  0                                                                  \\
        \sum_{k=0}^n (-1)^k \binom{n}{k} \binom{z}{k} \bigg/ \binom{y}{k}  
&  = &  \binom{y-z}{n} \bigg/ \binom{y}{n}                                 \\
        \sum_{k=0}^n (-1)^k \binom{n}{k} \binom{n+k}{k+1}                  
&  = &  0                                                                  \\*
        \sum_{k=0}^n (-1)^k \binom{n+k}{k} \binom{n+1}{k+1}                
&  = &  0                                                                  
\end{eqnarray*}
hold for every   $ n\in\N $   and all   $ y,z\in\R $.  
\end{lemma}

\bigskip
\begin{proof}
The first three identities can be found in 
Quaintance [2010; Equations (6.12), (6.44), (10.18)] 
and the last identity follows from the third. 
\end{proof} 

\bigskip
\begin{lemma}{} 
\label{l.app-1}
The identity 
$$  \int_{[0,u_{n+1}]} \dots \int_{[0,u_2]} S_n(\uuu) \,d\leb(u_1) \dots d\leb(u_n)  =  \frac{u_{n+1}^{2n}}{n!\,(n+1)!}  $$
holds for every   $ n\in\N $   and every   $ u_{n+1}\in\I $.  
\end{lemma}

\bigskip
\begin{proof}
We have 
\begin{eqnarray*}
        \int_{[0,u_2]} S_1(u_1) \,d\leb(u_1)  
&  = &  \int_{[0,u_2]}     u_1  \,d\leb(u_1)  
\;\,=\,\;    \frac{u_2^2}{2}                  
\end{eqnarray*}
which means that the assertion holds for   $ n=1 $. 
\\
Assume now that the assertion holds for all   $ k\in\{1,\dots,n\} $. 
Then we have 
\begin{eqnarray*}\lefteqn{
        \int_{[0,u_{n+2}]} \dots \int_{[0,u_2]} S_{n+1}(\uuu)                                                                            
        \,d\leb(u_1) \dots d\leb(u_{n+1})                                                                                               }\\*
&  = &  \int_{[0,u_{n+2}]} \dots \int_{[0,u_2]} \sum_{i=1}^{n+1} (-1)^{i+n+1} \frac{u_i^{n-i+2}}{(n-i+2)!}\, S_{i-1}(u_1,\dots,u_{i-1})  
        \,d\leb(u_1) \dots d\leb(u_{n+1})                                                                                                \\
&  = &  \sum_{i=1}^{n+1} (-1)^{i+n+1}                                                                                                    
        \int_{[0,u_{n+2}]} \dots \int_{[0,u_{i+1}]} \frac{u_i^{n-i+2}}{(n-i+2)!}                                                         \\*
&&                                                                                                                                       
        \Biggl[ \int_{[0,u_i]} \dots \int_{[0,u_2]} S_{i-1}(u_1,...,u_{i-1}) \,d\leb(u_1) \dots d\leb(u_{i-1}) \Biggr]                   
        d\leb(u_i) \dots d\leb(u_{n+1})                                                                                                  \\
&  = &  \sum_{i=1}^{n+1} (-1)^{i+n+1}                                                                                                    
        \int_{[0,u_{n+2}]} \dots \int_{[0,u_{i+1}]} \frac{u_i^{n-i+2}}{(n-i+2)!}\,\frac{u_i^{2(i-1)}}{(i-1)!\,i!}                        
        \,d\leb(u_i) \dots d\leb(u_{n+1})                                                                                                \\
&  = &  \sum_{i=1}^{n+1} (-1)^{i+n+1}                                                                                                    
        \int_{[0,u_{n+2}]} \dots \int_{[0,u_{i+1}]} \frac{u_i^{n+i}}{(n-i+2)!\,(i-1)!\,i!}                                               
        \,d\leb(u_i) \dots d\leb(u_{n+1})                                                                                                \\
&  = &  \sum_{i=1}^{n+1} (-1)^{i+n+1} \frac{u_{n+2}^{2n+2}}{(n-i+2)!\,(i-1)!\,i!}\,\frac{(n+i)!}{(2n+2)!}                                \\
&  = &  \frac{u_{n+2}^{2n+2}}{(n+1)!\,(n+2)!}                                                                                            
        \sum_{k=0}^n (-1)^{k+n+2} \frac{(n+1)!\,(n+2)!}{(n-k+1)!\,k!\,(k+1)!}\,\frac{(n+k+1)!}{(2n+2)!}                                  \\
&  = &  \frac{u_{n+2}^{2n+2}}{(n+1)!\,(n+2)!}                                                                                            
        \sum_{k=0}^n (-1)^k \binom{n+1}{k} \binom{n+1+k}{k+1} \bigg/ \biggl( (-1)^{n+2} \binom{2n+2}{n+2}  \biggr)                       \\*
&  = &  \frac{u_{n+2}^{2n+2}}{(n+1)!\,(n+2)!}                                                                                            
\end{eqnarray*}
In the final step we have used the third identity of Lemma 
\ref{l.combi}. 
\end{proof}

\bigskip
For   $ n\in\N $   and   $ m\in\{1,\dots,n\} $   we define the map   $ S_{n,m} : \I^n\to\R $   by letting 
\begin{eqnarray*}
        S_{n,m}(\uuu)                     
& := &  S_n(u_1,\dots,u_{n-m},1,\dots,1)  
\end{eqnarray*}
We have the following result: 
\pagebreak

\bigskip
\begin{lemma}{}
\label{l.app-2}
The identity 
\begin{eqnarray*} 
        \lefteqn{\int_{[0,u_{n+1}]} \dots \int_{[0,u_2]} S_{n,m}(\uuu) \,d\leb(u_1) \dots d\leb(u_n)}                                            \\*
&  = &  \sum_{l=0}^m \frac{u_{n+1}^{2n-l}}{(2n)!}\,\frac{1}{n-l+1}\,\binom{2n}{l} \sum_{h=0}^{m-l} (-1)^h \binom{2n-2l-h}{n-l} \binom{n-l+1}{h}  
\end{eqnarray*}
holds for 
every   $ n\in\N $,   
every   $ m\in\{1,\dots,n\} $   and 
every   $ u_{n+1}\in\I $. 
\end{lemma}

\bigskip
\begin{proof}
We have 
\begin{eqnarray*} 
        \int_{[0,u_2]} S_{1,1}(u_1) \,d\leb(u_1)                       
&  = &  \int_{[0,u_2]} S_1    (1)   \,d\leb(u_1)                       \\
&  = &  \int_{[0,u_2]}              \,d\leb(u_1)                       \\[2ex]
&  = &  u_2                                                            \\*[1ex]
&  = &  \sum_{l=0}^1 \frac{u_2^{2-l}}{2}\,\frac{1}{2-l}\,\binom{2}{l}  
        \sum_{h=0}^{1-l} (-1)^h \binom{2-2l-h}{1-l} \binom{2-l}{h}     
\end{eqnarray*}
which means that the identity holds for   $ n=1 $   and   $ m=1 $   and every   $ u_2\in\I $. 
\\
Assume henceforth that   $ n\geq2 $. 
\\
Let us first consider the case   $ m=1 $. 
In this case we have 
\begin{eqnarray*}
\lefteqn{                                                                                                                                 
        \int_{[0,u_{n+1}]} \dots \int_{[0,u_2]}                               S_{n,1}(\uuu)                \,d\leb(u_1) \dots d\leb(u_n) }\\*
&  = &  \int_{[0,u_{n+1}]} \dots \int_{[0,u_2]}                               S_n    (u_1,\dots,u_{n-1},1) \,d\leb(u_1) \dots d\leb(u_n)  \\
&  = &  \sum_{i=1}^{n-1} (-1)^{i+n}                                                                                                       
        \int_{[0,u_{n+1}]} \dots \int_{[0,u_2]} \frac{u_i^{n-i+1}}{(n-i+1)!}\,S_{i-1}(u_1,\dots,u_{i-1}  ) \,d\leb(u_1) \dots d\leb(u_n)  \\*
&& +    \int_{[0,u_{n+1}]} \dots \int_{[0,u_2]}                               S_{n-1}(u_1,\dots,u_{n-1}  ) \,d\leb(u_1) \dots d\leb(u_n)  
\end{eqnarray*}
For every   $ i\in\{1,\dots,n-1\} $, 
Lemma 
\ref{l.app-1} yields 
\begin{eqnarray*}
\lefteqn{                                                                                                                                           
        \int_{[0,u_{n+1}]} \dots \int_{[0,u_2]}     \frac{u_i^{n-i+1}}{(n-i+1)!}\,S_{i-1}(u_1,\dots,u_{i-1})     \,d\leb(u_1) \dots d\leb(u_n)     }\\
&  = &  \int_{[0,u_{n+1}]} \dots \int_{[0,u_{i+1}]} \frac{u_i^{n-i+1}}{(n-i+1)!}                                                                    \\*
&    &  \Biggl[                                                                                                                                     
        \int_{[0,u_i    ]} \dots \int_{[0,u_2]}                                   S_{i-1}(u_1,\dots,u_{i-1})     \,d\leb(u_1) \dots d\leb(u_{i-1})  
        \Biggr]                                                                                                    d\leb(u_i) \dots d\leb(u_n)      \\
&  = &  \int_{[0,u_{n+1}]} \dots \int_{[0,u_{i+1}]} \frac{u_i^{n-i+1}}{(n-i+1)!} \frac{u_i^{2(i-1)}}{(i-1)!\,i!} \,d\leb(u_i) \dots d\leb(u_n)      \\
&  = &  \frac{1}{(n-i+1)!(i-1)!\,i!} \int_{[0,u_{n+1}]} \dots \int_{[0,u_{i+1}]} u_i^{n+i-1}                     \,d\leb(u_i) \dots d\leb(u_n)      \\*
&  = &  \frac{1}{(n-i+1)!(i-1)!\,i!} \,\frac{(n+i-1)!}{(2n)!} \,u_{n+1}^{2n}                                                                        
\end{eqnarray*}
and summation yields 
\begin{eqnarray*}
\lefteqn{                                                                                                                                 
        \sum_{i=1}^{n-1} (-1)^{i+n}                                                                                                       
        \int_{[0,u_{n+1}]} \dots \int_{[0,u_2]} \frac{u_i^{n-i+1}}{(n-i+1)!}\,S_{i-1}(u_1,\dots,u_{i-1}  ) \,d\leb(u_1) \dots d\leb(u_n) }\\*
&  = &  \sum_{i=1}^{n-1} (-1)^{i+n} \frac{1}{(n-i+1)!(i-1)!\,i!} \,\frac{(n+i-1)!}{(2n)!} \,u_{n+1}^{2n}                                  \\
&  = &  (-1)^{n+1}\,\frac{u_{n+1}^{2n}}{(2n)!   }         \sum_{k=0  }^{n-2} (-1)^k    \,\frac{(n+k)!}{(n-k)!k!\,(k+1)!}                  \\
&  = &  (-1)^{n+1}\,\frac{u_{n+1}^{2n}}{(2n)!\,n}         \sum_{k=0  }^{n-2} (-1)^k    \,\binom{n}{k} \binom{n+k}{k+1}                    \\
&  = &  (-1)^{n+1}\,\frac{u_{n+1}^{2n}}{(2n)!\,n}    (-1) \sum_{k=n-1}^n     (-1)^k    \,\binom{n}{k} \binom{n+k}{k+1}                    \\*
&  = &  \frac{u_{n+1}^{2n}}{n!\,(n+1)!} - \frac{u_{n+1}^{2n}}{(n-1)!\,n!\,2n}                                                             
\end{eqnarray*}
We also obtain 
\begin{eqnarray*}
\lefteqn{ 
        \int_{[0,u_{n+1}]} \dots \int_{[0,u_2]} S_{n-1}(u_1,\dots,u_{n-1}) \,d\leb(u_1) \dots d\leb(u_n)                        }\\
&  = &  \int_{[0,u_{n+1}]}                                                                                                       
        \Biggl[                                                                                                                  
        \int_{[0,u_n    ]} \dots \int_{[0,u_2]} S_{n-1}(u_1,\dots,u_{n-1}) \,d\leb(u_1) \dots d\leb(u_{n-1}) \Biggr] d\leb(u_n)  \\
&  = &  \int_{[0,u_{n+1}]}                      \frac{u_n^{2(n-1)}}{(n-1)!\,n!}                                    \,d\leb(u_n)  \\*
&  = &  \frac{u_{n+1}^{2n-1}}{(n-1)!\,n!\,(2n-1)}                                                                                
\end{eqnarray*}
Now summation yields 
\begin{eqnarray*}
\lefteqn{                                                                                                                                        
        \int_{[0,u_{n+1}]} \dots \int_{[0,u_2]}                               S_{n,1}(\uuu             )          \,d\leb(u_1) \dots d\leb(u_n) }\\
&  = &  \sum_{i=1}^{n-1} (-1)^{i+n}                                                                                                              
        \int_{[0,u_{n+1}]} \dots \int_{[0,u_2]} \frac{u_i^{n-i+1}}{(n-i+1)!}\,S_{i-1}(u_1,\dots,u_{i-1})          \,d\leb(u_1) \dots d\leb(u_n)  \\
&& +    \int_{[0,u_{n+1}]} \dots \int_{[0,u_2]}                               S_{n-1}(u_1,\dots,u_{n-1})          \,d\leb(u_1) \dots d\leb(u_n)  \\
&  = &  \frac{u_{n+1}^{2n}}{n!\,(n+1)!} - \frac{u_{n+1}^{2n}}{(n-1)!\,n!\,2n} + \frac{u_{n+1}^{2n-1}}{(2n-1)(n-1)!\,n!}                          \\*
&  = &  \sum_{l=0}^1 \frac{u_{n+1}^{2n-l}}{(2n)!}\,\frac{1}{n-l+1}\,\binom{2n}{l} \sum_{h=0}^{1-l} (-1)^h \binom{2n-2l-h}{n-l} \binom{n-l+1}{h}  
\end{eqnarray*}
Let us now consider the case   $ m\in\{2,\dots,n\} $   and assume that the identity 
\begin{eqnarray*} 
        \lefteqn{\int_{[0,u_{q+1}]} \dots \int_{[0,u_2]} S_{q,p}(\uuu) \,d\leb(u_1) \dots d\leb(u_q)}                                            \\*
&  = &  \sum_{l=0}^p \frac{u_{q+1}^{2q-l}}{(2q)!}\,\frac{1}{q-l+1}\,\binom{2q}{l} \sum_{h=0}^{p-l} (-1)^h \binom{2q-2l-h}{q-l} \binom{q-l+1}{h}  
\end{eqnarray*}
holds for 
every   $ q\in\{1,\dots,n-1\} $,   
every   $ p\in\{1,\dots,q\} $   and 
every   $ u_{q+1}\in\I $.  
We have 
\begin{eqnarray*}
\lefteqn{                                                                                                                                
        \int_{[0,u_{n+1}]} \dots \int_{[0,u_2]} S_{n,m}(\uuu)                                            \, d\leb(u_1) \dots d\leb(u_n) }\\*
&  = &  \sum_{i=1}^{n-m} (-1)^{i+n} \int_{[0,u_{n+1}]} \dots \int_{[0,u_2]} \frac{u_i^{n-i+1}}{(n-i+1)!} \,                              
                                                                               S_{i-1}(u_1,\dots,u_{i-1})\, d\leb(u_1) \dots d\leb(u_n)  \\
&     & +\,(-1)^{2n-m+1} \int_{[0,u_{n+1}]} \dots \int_{[0,u_2]} \frac{1}{m!}\,S_{n-m}(u_1,\dots,u_{n-m})\, d\leb(u_1) \dots d\leb(u_n)  \\*
&     & +\, \sum_{i=n-m+2}^{n} (-1)^{i+n} \int_{[0,u_{n+1}]} \dots \int_{[0,u_2]} \frac{1}{(n-i+1)!}\,                                   
                                                                       S_{i-1,m-n+i-1}(u_1,\dots,u_{i-1})\, d\leb(u_1) \dots d\leb(u_n)  
\end{eqnarray*}
For every $ i \in \{1,...,n-m\} $, 
Lemma 
\ref{l.app-1} yields 
\begin{eqnarray*}
\lefteqn{                                                                                                                                                   
        \int_{[0,u_{n+1}]} \dots \int_{[0,u_2]}     \frac{u_i^{n-i+1}}{(n-i+1)!}\,S_{i-1}(u_1,\dots,u_{i-1})     \,d\leb(u_1) \dots d\leb(u_n    )         }\\*
&  = &  \int_{[0,u_{n+1}]} \dots \int_{[0,u_{i+1}]} \frac{u_i^{n-i+1}}{(n-i+1)!}                                                                            \\*
&&      \Biggl[ \int_{[0,u_i]} \dots \int_{[0,u_2]}                               S_{i-1}(u_1,\dots,u_{i-1})     \,d\leb(u_1) \dots d\leb(u_{i-1}) \Biggr]  
                                                                                                                   d\leb(u_i) \dots d\leb(u_n    )          \\
&  = &  \int_{[0,u_{n+1}]} \dots \int_{[0,u_{i+1}]} \frac{u_i^{n-i+1}}{(n-i+1)!}\,\frac{u_i^{2(i-1)}}{(i-1)!\,i!}\,d\leb(u_i) \dots d\leb(u_n    )          \\ 
&  = &  \int_{[0,u_{n+1}]} \dots \int_{[0,u_{i+1}]} \frac{u_i^{n+i-1}}{(n-i+1)!(i-1)!\,i!}                       \,d\leb(u_i) \dots d\leb(u_n    )          \\*
&  = &  \frac{u_{n+1}^{2n}}{(n-i+1)!\,(i-1)!\,i!}\,\frac{(n+i-1)!}{(2n)!}                                                                                   
\end{eqnarray*}
and summation yields 
\begin{eqnarray*}
\lefteqn{                                                                                                                              
        \sum_{i=1}^{n-m} (-1)^{i+n} \int_{[0,u_{n+1}]} \dots \int_{[0,u_2]} \frac{u_i^{n-i+1}}{(n-i+1)!}\, S_{i-1}(u_1,\dots,u_{i-1})  
        \,d\leb(u_1) \dots d\leb(u_n)                                                                                                 }\\
&  = &  \sum_{i=1}^{n-m} (-1)^{i+n} \frac{u_{n+1}^{2n}}{(n-i+1)!\,(i-1)!\,i!}\,\frac{(n+i-1)!}{(2n)!}      \hspace*{40mm}                            \\
&  = &  \frac{u_{n+1}^{2n}}{(2n)!}\,\frac{1}{n+1}\,(-1)^{n+1} \sum_{k=0}^{n-m-1} (-1)^k \binom{n+k }{k} \binom{n+1}{k+1}               \\
&  = &  \frac{u_{n+1}^{2n}}{(2n)!}\,\frac{1}{n+1}\,(-1)^{n+2} \sum_{k=n-m}^n     (-1)^k \binom{n+k }{k} \binom{n+1}{k+1}               \\
&  = &  \frac{u_{n+1}^{2n}}{(2n)!}\,\frac{1}{n+1}             \sum_{h=0}^m       (-1)^h \binom{2n-h}{n} \binom{n+1}{h}                 
\end{eqnarray*}
Because of Lemma 
\ref{l.app-1} we also obtain 
\begin{eqnarray*}
\lefteqn{                                                                                                                             
        (-1)^{2n-m+1} \int_{[0,u_{n+1}]} \dots \int_{[0,u_2]} \frac{1}{m!}\,S_{n-m}(u_1,\dots,u_{n-m})                                
                                                                                     \,d\leb(u_1      ) \dots d\leb(u_n    )         }\\*
&  = &  (-1)^{m-1}\,\frac{1}{m!} \int_{[0,u_{n+1}]} \dots \int_{[0,u_{n-m+2}]}                                                        \\*
&&      \Biggl[ \int_{[0,u_{n-m+1}]} \dots \int_{[0,u_2]} S_{n-m}(u_1,\dots,u_{n-m}) \,d\leb(u_1      ) \dots d\leb(u_{n-m}) \Biggr]  
                                                                                       d\leb(u_{n-m+1}) \dots d\leb(u_n    )          \\
&  = &  (-1)^{m-1}\,\frac{1}{m!} \int_{[0,u_{n+1}]} \dots \int_{[0,u_{n-m+2}]} \frac{u_{n-m+1}^{2(n-m)}}{(n-m)!\,(n-m+1)!}\,          
                                                                                       d\leb(u_{n-m+1}) \dots d\leb(u_n    )          \\
&  = &  (-1)^{m-1}\,\frac{u_{n+1}^{2n-m}}{m!\,(n-m)!\,(n-m+1)!}\,\frac{(2n-2m)!}{(2n-m)!}                     \\*
&  = &  (-1)^{m-1}\,\frac{u_{n+1}^{2n-m}}{(2n)!}\,\frac{1}{n-m+1} \binom{2n}{m} \binom{2n-2m}{n-m}
\end{eqnarray*}
Furthermore, 
we obtain 
\begin{eqnarray*}
\lefteqn{
        \sum_{i=n-m+2}^n (-1)^{i+n}                                                                                               
        \int_{[0,u_{n+1}]} \dots \int_{[0,u_2]} \frac{1}{(n-i+1)!}\,                                                              
        S_{i-1,m-n+i-1}(u_1,\dots,u_{i-1})                                             \,d\leb(u_1) \dots d\leb(u_n)             }\\
&  = &  \sum_{i=n-m+2}^n (-1)^{i+n} \frac{1}{(n-i+1)!}                                                                            
        \int_{[0,u_{n+1}]} \dots \int_{[0,u_{i+1}]}                                                                               \\*
&&      \Biggl[ \int_{[0,u_i]} \dots \int_{[0,u_2]} S_{i-1,m-n+i-1}(u_1,\dots,u_{i-1}) \,d\leb(u_1) \dots d\leb(u_{i-1}) \Biggr]  
                                                                                         d\leb(u_i) \dots d\leb(u_n)              \\
&  = &  \sum_{i=n-m+2}^n (-1)^{i+n} \frac{1}{(n-i+1)!}                                                                            
        \int_{[0,u_{n+1}]} \dots \int_{[0,u_{i+1}]}                                                                               \\*
&&      \Biggl[ \sum_{l=0}^{m-n+i-1} \frac{u_{i}^{2i-2-l}}{(2i-2)!}\,\frac{1}{i-l} \binom{2i-2}{l}                                
        \sum_{h=0}^{m-n+i-1-l} (-1)^h \binom{2i-2-2l-h}{i-1-l} \binom{i-l}{h} \Biggr]    d\leb(u_i) \dots d\leb(u_n)              \\
&  = &  \sum_{i=n-m+2}^n \sum_{l=0}^{m-n+i-1} \sum_{h=0}^{m-n+i-1-l} (-1)^{i+n+h} \frac{1}{i-l}\,                                 
        \frac{1}{(n-i+1)!\,(2i-2-l)!\,l!}                                                                                         \\*
&&      \binom{2i-2-2l-h}{i-1-l} \binom{i-l}{h}                                                                                   
        \int_{[0,u_{n+1}]} \dots \int_{[0,u_{i+1}]} u_{i}^{2i-2-l}                     \,d\leb(u_i) \dots d\leb(u_n)              \\
&  = &  \sum_{i=n-m+2}^n \sum_{l=0}^{m-n+i-1} \sum_{h=0}^{m-n+i-1-l} (-1)^{i+n+h}                                                 
        \frac{1}{(n-i+1)!\,(2i-2-l)!\,l!}                                                                                         \\*
&&      \frac{(2i-2-2l-h)!}{(i-1-l-h)!\,(i-l-h)!\,h!}\,\frac{(2i-2-l)!}{(n+i-1-l)!}\,u_{n+1}^{n+i-1-l}                            \\
&  = &  \sum_{k=1}^{m-1} \sum_{l=0}^k   \sum_{h=0}^{k-l} (-1)^{k-m+1+h} \frac{1}{(m-k)!\, l!}                                     \\*
&&      \frac{(2n-2m+2k-2l-h)!}{(n-m+k-l-h)!\,(n-m+k+1-l-h)!\,h!}\,\frac{1}{(2n-m+k-l)!}\,u_{n+1}^{2n-m+k-l}                      \\
&  = &  \sum_{k=1}^{m-1} \sum_{j=m-k}^m \sum_{h=0}^{m-j} (-1)^{k-m+1+h} \frac{1}{(m-k)!\,(k-m+j)!}                                \\*
&&      \frac{(2n-2j-h)!}{(n-j-h)!\,(n-j+1-h)!\,h!}\,\frac{1}{(2n-j)!}\,u_{n+1}^{2n-j}                                            \\
&  = &  \sum_{k=1}^{m-1} \sum_{j=m-k}^m \frac{u_{n+1}^{2n-j}}{(2n)!} \binom{2n}{j}                                                
        \sum_{h=0}^{m-j} (-1)^{k-m+1+h} \frac{j!}{(m-k)!\,(k-m+j)!}\,\frac{(2n-2j-h)!}{(n-j-h)!\,(n-j+1-h)!\,h!}                  \\
&  = &  \sum_{j=1}^{m-1} \sum_{k=m-j}^{m-1} \frac{u_{n+1}^{2n-j}}{(2n)!} \binom{2n}{j}                                            
        \sum_{h=0}^{m-j} (-1)^{k-m+1+h} \binom{j}{m-k} \binom{2n-2j-h}{n-j} \binom{n-j+1}{h} \frac{1}{n-j+1}                      \\*
&&   +\,(-1)^{m+1} \, \frac{u_{n+1}^{2n-m}}{(2n)!}\,\frac{1}{n-m+1} \binom{2n}{m} \binom{2n-2m}{n-m}                              
        \sum_{k=1}^{m-1} (-1)^k         \binom{m}{k}                                                                              \\
&  = &  \sum_{j=1}^{m-1} \frac{u_{n+1}^{2n-j}}{(2n)!}\,\frac{1}{n-j+1} \binom{2n}{j}                                              
        \sum_{h=0}^{m-j} (-1)^h         \binom{2n-2j-h}{n-j} \binom{n-j+1}{h}                                                     
        \sum_{k=m-j}^{m-1} (-1)^{k-m+1} \binom{j}{m-k}                                                                            \\*
&&   +\,(-1)^{m+1} \,\frac{u_{n+1}^{2n-m}}{(2n)!}\,\frac{1}{n-m+1} \binom{2n}{m} \binom{2n-2m}{n-m}                               
        \sum_{k=1}^{m-1} (-1)^k         \binom{m}{k}                                                                              \\
&  = &  \sum_{j=1}^{m-1} \frac{u_{n+1}^{2n-j}}{(2n)!}\,\frac{1}{n-j+1} \binom{2n}{j}                                              
        \sum_{h=0}^{m-j} (-1)^h         \binom{2n-2j-h}{n-j} \binom{n-j+1}{h}                                                     \\*
&&   +\,(-1)^m \, \frac{u_{n+1}^{2n-m}}{(2n)!}\,\frac{1}{n-m+1} \binom{2n}{m} \binom{2n-2m}{n-m}                                  
                + \frac{u_{n+1}^{2n-m}}{(2n)!}\,\frac{1}{n-m+1} \binom{2n}{m} \binom{2n-2m}{n-m}                                  
\end{eqnarray*}
Combining the previous results we obtain 
\begin{eqnarray*}
\lefteqn{
        \int_{[0,u_{n+1}]} \dots \int_{[0,u_2]} S_{n,m}(\uuu)            \,d\leb(u_1) \dots d\leb(u_n)   }\\
&  = &  \sum_{i=1}^{n-m} (-1)^{i+n} \int_{[0,u_{n+1}]} \dots \int_{[0,u_2]}                               
        \frac{u_i^{n-i+1}}{(n-i+1)!}\,S_{i-1}(u_1,\dots,u_{i-1})         \,d\leb(u_1) \dots d\leb(u_{n})  \\*
&&   +\,(-1)^{2n-m+1} \int_{[0,u_{n+1}]} \dots \int_{[0,u_2]}                                             
        \frac{1}{m!}                \,S_{n-m}(u_1,\dots,u_{n-m})         \,d\leb(u_1) \dots d\leb(u_n)    \\*
&&   +\,\sum_{i=n-m+2}^{n} (-1)^{i+n} \int_{[0,u_{n+1}]} \dots \int_{[0,u_2]}                             
        \frac{1}{(n-i+1)!}          \,S_{i-1,m-n+i-1}(u_1,\dots,u_{i-1}) \,d\leb(u_1) \dots d\leb(u_n)    \\
&  = &  \frac{u_{n+1}^{2n}}{(2n)!}\,\frac{1}{n+1} \sum_{h=0}^m (-1)^h  \binom{2n-h}{n} \binom{n+1}{h}     \\*
&&   +\,(-1)^{m-1}     \,\frac{u_{n+1}^{2n-m}}{(2n)!}\,\frac{1}{n-m+1} \binom{2n}{m} \binom{2n-2m}{n-m}   \\*
&&   +\,\sum_{j=1}^{m-1} \frac{u_{n+1}^{2n-j}}{(2n)!}\,\frac{1}{n-j+1} \binom{2n}{j}                      
        \sum_{h=0}^{m-j} (-1)^h \binom{2n-2j-h}{n-j} \binom{n-j+1}{h}                                     \\*
&&   +\,(-1)^m         \,\frac{u_{n+1}^{2n-m}}{(2n)!}\,\frac{1}{n-m+1} \binom{2n}{m}\binom{2n-2m}{n-m}    
                       + \frac{u_{n+1}^{2n-m}}{(2n)!}\,\frac{1}{n-m+1} \binom{2n}{m}\binom{2n-2m}{n-m}    \\*
&  = &  \sum_{l=0}^m     \frac{u_{n+1}^{2n-l}}{(2n)!}\,\frac{1}{n-l+1} \binom{2n}{l}                      
        \sum_{h=0}^{m-l} (-1)^h \binom{2n-2l-h}{n-l} \binom{n-l+1}{h}                                     
\end{eqnarray*}
This proves the identity for   $ n+1 $   and   $ m\in\{2,\dots,n+1\} $. 
\end{proof}

\bigskip
\begin{corollary}{} 
\label{c.app-1}
The identity 
\begin{eqnarray*}
\lefteqn{                                                                                                      
        \int_{[0,1]} \dots \int_{[0,u_2]} S_{n,m}(\uuu) \,d\leb(u_1) \dots d\leb(u_n)                         }\\
&  = &  \frac{1}{n!\,(n+1)!} + \frac{1}{(2n)!} \sum_{k=1}^m \frac{(2n-2k)!}{(n-k)!\,(n+1-k)!} \binom{2k-1}{k}  
\end{eqnarray*}
holds for every   $ n\in\N $   and every   $ m\in\{1,\dots,n\} $. 
\end{corollary}

\bigskip
\begin{proof}
From Lemmas 
\ref{l.app-2} and 
\ref{l.combi} we obtain
\begin{eqnarray*} 
\lefteqn{                                                                                                                   
        \int_{[0,1]} \dots \int_{[0,u_2]} S_{n,m}(\uuu) \,d\leb(u_1) \dots d\leb(u_n)                                      }\\*
&  = &  \sum_{l=0}^m \frac{1}{(2n)!}\,\frac{1}{n-l+1}\,\binom{2n}{l}                                                        
        \sum_{h=0}^{m-l} (-1)^h \binom{2n-2l-h}{n-l} \binom{n-l+1}{h}                                                       \\
&  = &  \sum_{l=0}^m            \frac{1}{(2n)!}\,\frac{1}{n-l+1}\,\binom{2n}{l}                                             
        \sum_{k=l}^m (-1)^{k-l} \binom{2n-l-k }{n-l} \binom{n-l+1}{k-l}                                                     \\
&  = &  \sum_{k=0}^m (-1)^k     \frac{1}{(2n)!\,(n-k)!\,(n+1-k)!}                                                           
        \sum_{l=0}^k (-1)^l     \binom{2n}{l}\frac{(2n-l-k)!}{(k-l)!}                                                       \\
&  = &  \sum_{k=0}^m (-1)^k     \frac{(2n-k)!}{(2n)!\,(n-k)!\,(n+1-k)!\,k!}                                                 
        \sum_{l=0}^k (-1)^l     \binom{k}{l}\binom{2n}{l} \bigg/ \binom{2n-k}{l}                                            \\
&  = &  \sum_{k=0}^m (-1)^k     \frac{(2n-k)!}{(2n)!\,(n-k)!\,(n+1-k)!\,k!}\,(-1)^k \binom{2k-1}{k} \bigg/ \binom{2n-k}{k}  \\*
&  = &  \frac{1}{n!\,(n+1)!} + \frac{1}{(2n)!} \sum_{k=1}^{m} \frac{(2n-2k)!}{(n-k)!\,(n+1-k)!}            \binom{2k-1}{k}  
\end{eqnarray*}
This proves the assertion.
\end{proof}


\section*{Acknowledgement}

The authors would like to thank Georg Berschneider and Klaus Th.\ Hess for drawing their attention to 
Cahill et al.\ [2002] and the collection of combinatorial formulas by 
Quaintance [2010].


\section*{References}
\small 

Av{\'e}rous J, Genest C, Kochar SC [2005]: 
On the dependence structure of order statistics. 
Journal of Multivariate Analysis 94, 159--171. 

\smallskip
Cahill ND, D'Errico JR, Narayan DA, Narayan JY [2002]: 
Fibonacci Determinants. 
The College Mathematical Journal 33, 221--225

\smallskip
Cap{\'e}ra{\`a} P, Foug{\`e}res AL, Genest C [1997]: 
A stochastic ordering based on a decomposition of Kendall's tau. 
In: 
Bene\u{s} V, \u{S}tep{\'a}n J (eds.), 
Distribution with Given Marginals and Moment Problems, 
pp.\ 81--86.
Dordrecht, Kluwer Academic Publishers. 

\smallskip
Dietz M, Fuchs S, Schmidt KD [2016]: 
On order statistics and their copulas. 
Statistics and Probability Letters 117, 165--172. 

\smallskip
Durante F, Sempi C [2016]: 
Principles of Copula Theory.
London, Chapman \& Hall. 

\smallskip
Fuchs S [2014]: 
Multivariate copulas: Transformations, symmetry, order and measures of concordance.
Kybernetika 50, 725--743. 

\smallskip
Fuchs S [2016]: 
A biconvex form for copulas.
Dependence Modeling 4, 63--75. 

\smallskip
Fuchs S, McCord Y, Schmidt KD [2018]: 
Characterizations of copulas attaining the bounds of multivariate Kendall's tau. 
Journal of Optimization Theory and Applications. 

\smallskip
Navarro J, Balakrishnan N [2010]: 
Study of some measures of dependence between order statistics and systems. 
Journal of Multivariate Analysis 101, 52--67. 

\smallskip
Navarro J, Spizzichino F [2010]: 
On the relationship between copulas of order statistics and marginal distributions. 
Statistics and Probability Letters 80, 473--479. 

\smallskip
Nelsen RB [2002]: 
Concordance and copulas: A survey.
In: Cuadras CM, Fortiana J, Rodriguez--Lallena JA (eds.), 
Distributions with Given Marginals and Statistical Modelling, pp.\ 169--177. 
Dordrecht, Kluwer Academic Publishers. 

\smallskip
Nelsen RB [2006]: 
An Introduction to Copulas. Second edition. 
New York, Springer. 

\smallskip
Quaintance J [2010]: 
Combinatorial Identities: Table I: Intermediate Techniques for Summing Finite Series. 
https://www.math.wvu.edu/~gould/Vol.4.PDF

\bigskip
\vfill\hspace*{\fill}\today

\end{document}